\documentclass[a4paper,11pt]{amsart}
\usepackage{import}
\calclayout

\usepackage[utf8]{inputenc}
\usepackage{enumitem}
\usepackage[symbol]{footmisc}
\usepackage{graphicx}
\usepackage{xcolor}
\usepackage{float}
\usepackage{tikz-cd}
\usepackage{tikz}
\usepackage{adjustbox}
\usepackage{standalone}
\usepackage{relsize}
\usepackage{xfrac}
\usepackage[nospace,noadjust]{cite}
\usepackage{amsthm}
\usepackage{upgreek}
\usepackage{bm}
\usepackage{comment}
\usepackage[T1]{fontenc} 

\usetikzlibrary{arrows.meta,patterns}
\usetikzlibrary{ipe} 

\setenumerate[1]{label={\roman*)}}

\theoremstyle{plain}
\newtheorem{Theorem}{Theorem}[section]
\newtheorem{Proposition}[Theorem]{Proposition}
\newtheorem*{Theorem*}{Theorem}
\newtheorem{Corollary}[Theorem]{Corollary}
\newtheorem{Lemma}[Theorem]{Lemma}

\newtheorem{innercustomthm}{Theorem}
\newenvironment{customthm}[1]
{\renewcommand\theinnercustomthm{#1}\innercustomthm}
{\endinnercustomthm}

\theoremstyle{remark}
\newtheorem{Remark}[Theorem]{Remark}

\newtheorem{Example}[Theorem]{Example}

\theoremstyle{definition}
\newtheorem{Definition}[Theorem]{Definition}

\renewcommand{\S}{\ensuremath{\mathbb{S}}}
\newcommand{\R}{\ensuremath{\mathbb{R}}}
\newcommand{\Z}{\ensuremath{\mathbb{Z}}}
\newcommand{\C}{\ensuremath{\mathbb{C}}}
\renewcommand{\H}{\ensuremath{\mathbb{H}}}
\newcommand{\V}{\ensuremath{\mathbb{V}}}

\newcommand{\lspan}{\bm{\langle}}
\newcommand{\rspan}{\bm{\rangle}}
\newcommand{\sdfrac}[2]{\mbox{\small$\displaystyle\frac{#1}{#2}$}}
\newcommand\Wtilde{\stackrel{\sim}{\smash{\mathcal{W}}\rule{0pt}{1.1ex}}}
\newcommand\Dtilde{\stackrel{\sim}{\smash{\mathcal{D}}\rule{0pt}{1.1ex}}}
\newcommand\Etilde{\stackrel{\sim}{\smash{\mathcal{E}}\rule{0pt}{1.1ex}}}

\newcommand{\orlspan}{\bm{\textlangle}}
\newcommand{\orrspan}{\bm{\textrangle}}

\DeclareCollectionInstance{plainmath}{xfrac}{mathdefault}{math}
{
 slash-symbol = \text{\relscale{1.2}$/$},
}
\UseCollection{xfrac}{plainmath}

\tikzstyle{ipe stylesheet} = [
  ipe import,
  even odd rule,
  line join=round,
  line cap=butt,
  ipe pen normal/.style={line width=0.4},
  ipe pen heavier/.style={line width=0.8},
  ipe pen fat/.style={line width=1.2},
  ipe pen ultrafat/.style={line width=2},
  ipe pen normal,
  ipe mark normal/.style={ipe mark scale=3},
  ipe mark large/.style={ipe mark scale=5},
  ipe mark small/.style={ipe mark scale=2},
  ipe mark tiny/.style={ipe mark scale=1.1},
  ipe mark normal,
  /pgf/arrow keys/.cd,
  ipe arrow normal/.style={scale=7},
  ipe arrow large/.style={scale=10},
  ipe arrow small/.style={scale=5},
  ipe arrow tiny/.style={scale=3},
  ipe arrow normal,
  /tikz/.cd,
  ipe arrows, 
  <->/.tip = ipe normal,
  ipe dash normal/.style={dash pattern=},
  ipe dash dotted/.style={dash pattern=on 1bp off 3bp},
  ipe dash dashed/.style={dash pattern=on 4bp off 4bp},
  ipe dash dash dotted/.style={dash pattern=on 4bp off 2bp on 1bp off 2bp},
  ipe dash dash dot dotted/.style={dash pattern=on 4bp off 2bp on 1bp off 2bp on 1bp off 2bp},
  ipe dash normal,
  ipe node/.append style={font=\normalsize},
  ipe stretch normal/.style={ipe node stretch=1},
  ipe stretch TikZ-normal/.style={ipe node stretch=1},
  ipe stretch normal,
  ipe opacity 10/.style={opacity=0.1},
  ipe opacity 30/.style={opacity=0.3},
  ipe opacity 50/.style={opacity=0.5},
  ipe opacity 75/.style={opacity=0.75},
  ipe opacity opaque/.style={opacity=1},
  ipe opacity opaque,
]

\definecolor{red}{rgb}{1,0,0}
\definecolor{blue}{rgb}{0,0,1}
\definecolor{green}{rgb}{0,1,0}
\definecolor{yellow}{rgb}{1,1,0}
\definecolor{orange}{rgb}{1,0.647,0}
\definecolor{gold}{rgb}{1,0.843,0}
\definecolor{purple}{rgb}{0.627,0.125,0.941}
\definecolor{gray}{rgb}{0.745,0.745,0.745}
\definecolor{brown}{rgb}{0.647,0.165,0.165}
\definecolor{navy}{rgb}{0,0,0.502}
\definecolor{pink}{rgb}{1,0.753,0.796}
\definecolor{seagreen}{rgb}{0.18,0.545,0.341}
\definecolor{turquoise}{rgb}{0.251,0.878,0.816}
\definecolor{violet}{rgb}{0.933,0.51,0.933}
\definecolor{darkblue}{rgb}{0,0,0.545}
\definecolor{darkcyan}{rgb}{0,0.545,0.545}
\definecolor{darkgray}{rgb}{0.663,0.663,0.663}
\definecolor{darkgreen}{rgb}{0,0.392,0}
\definecolor{darkmagenta}{rgb}{0.545,0,0.545}
\definecolor{darkorange}{rgb}{1,0.549,0}
\definecolor{darkred}{rgb}{0.545,0,0}
\definecolor{lightblue}{rgb}{0.678,0.847,0.902}
\definecolor{lightcyan}{rgb}{0.878,1,1}
\definecolor{lightgray}{rgb}{0.827,0.827,0.827}
\definecolor{lightgreen}{rgb}{0.565,0.933,0.565}
\definecolor{lightyellow}{rgb}{1,1,0.878}
\definecolor{black}{rgb}{0,0,0}
\definecolor{white}{rgb}{1,1,1}

\title[Contact space of null geodesics: compact, Engel, retrievability]{On the space of
null geodesics of a \\
spacetime: the compact case, Engel\\ geometry and~retrievability}
\author{Adrià Marín-Salvador and Roberto Rubio}

\thanks{The first author has been supported by the Spanish MEFP under the \textit{Beca de Colaboración} BDNS 512590 and by ``la Caixa'' Foundation (ID 100010434) under the fellowship LCF/BQ/EU21/11890122. The second author has been supported by the European Union’s Horizon 2020 research and innovation programme under the Marie Sklodowska-Curie grant agreement No 750885 GENERALIZED and by the Spanish State Research Agency under the grant PID2019-109339GA-C32}

\address{A. Marín-Salvador \\ 
University of Oxford, Oxford OX2 6GG, United Kingdom, and
Universitat Aut\`onoma de Barcelona, 08193 Barcelona, Spain}

\email{adria.marin@maths.ox.ac.uk}

\address{R. Rubio, 
Universitat Aut\`onoma de Barcelona, 08193 Barcelona, Spain}

\email{roberto.rubio@uab.es}

\begin{document}
\begin{abstract}
We compute the contact manifold of null geodesics of the family of spacetimes $\{(\S^2\times\S^1, g_\circ-\frac{d^2}{c^2}dt^2)\}_{d,c\in\mathbb{N}^+\text{ coprime}}$, with $g_\circ$ the round metric on $\S^2$ and $t$ the $\S^1$-coordinate. We find that these are the lens spaces $L(2c,1)$ together with the pushforward of the canonical contact structure on $ST\S^2\cong L(2,1)$ under the natural projection $L(2,1)\to L(2c,1)$. We extend this computation to $Z\times \S^1$ for $Z$ a Zoll manifold. On the other hand, motivated by these examples, we show how Engel geometry can be used to describe the manifold of null geodesics of a certain class of three-dimensional spacetimes, by considering the Cartan deprolongation of their Lorentz prolongation. We characterize the three-dimensional contact manifolds that are contactomorphic to the space of null geodesics of a spacetime. The characterization consists in the existence of an overlying Engel manifold with a certain foliation and, in this case, we also retrieve the spacetime.
\end{abstract}
\maketitle

\vspace{-.45cm}
\section{Introduction}
A spacetime is a Lorentzian manifold together with a choice of a global timelike vector field, that is, a vector field of negative length at all points. For a spacetime, its space of null geodesics $\mathcal{N}$ consists of the family of unparametrized geodesics with null tangent vectors at all points \cite{Low_1988,Low_1990}. When $\mathcal{N}$ is a manifold, it can be equipped with a canonical contact structure $\mathcal{H}$ (see \cite{Low_1993, Low_2001}).

The contact structure $\mathcal{H}$ has proved to be essential in the theory, yielding important results on causality, providing, for instance, obstructions to two events (points) being on the same non-spacelike curve \cite{Nat_rio_2004, Chernov_2009, Chernov_2010}.  The spaces of null geodesics and their contact structures were computed explicitly in some noncompact cases \cite{Bautista2008, Godoy2013} and the question of whether two spacetimes with diffeomorphic spaces of null geodesics must be diffeomorphic was addressed, under the name of reconstruction, in \cite{Bautista_2014}. Until very recently, the only explicit cases where $\mathcal{N}$ was known to be a manifold were globally hyperbolic spacetimes, which are diffeomorphic to $C\times \R$ for a Cauchy hypersurface $C$ \cite{Bernal2003} and for which  $\mathcal{N}\cong STC$ \cite{Low_1993}, and Zoll (or Zollfrei)  manifolds \cite{guillemin, guillemin-book}, whose null geodesics are all periodic. In the last years this subject has attracted more attention with the negative answer to Guillemin's conjecture that every Zoll $3$-dimensional spacetime is covered by $\S^2\times \R$ \cite{suhr}, and new classes of examples for which  $\mathcal{N}$ is a manifold \cite{hedicke-suhr, hedicke}. However, to the best of our knowledge, there are hardly any explicit calculations of spaces of null geodesics and their contact structures for compact spacetimes or results on the possibility of retrieving the spacetime from its space of null geodesics. 

Firstly, we consider the spacetimes $(\S^2\times\S^1, g_{c/d})$ for the family of metrics $g_c = g_\circ-\frac{d^2}{c^2}dt^2$ with $c,d\in\mathbb{N}^+$, where $g_\circ$ is the round metric on $\S^2$ and $t$ is the angle coordinate on $\S^1$. Its space of null geodesics $\mathcal{N}_{c/d}$ and its contact structure are described in terms of the lens spaces $L(2c,1)$ for $\gcd(c,d) = 1$, by using a quaternionic approach to the Hopf fibration and $ST\S^2$ that we develop in Section \ref{sec:quaternionic}. We prove,
\begin{customthm}{\ref{Thm: A}}\label{Thm: IntroA}
For any $c,d\geq 1$ with $\gcd(c,d)=1$,
\[
\mathcal{N}_{c/d}\cong L(2c,1).
\]
\end{customthm}

\begin{customthm}{\ref{Thm: B}}\label{Thm: IntroB}
The canonical contact structure $\mathcal{H}$ on $ST\S^2\cong\mathcal{N}_{1/d}$ is the canonical contact structure $\chi$ on $ST\S^2$. In general, the canonical contact structure $\mathcal{H}$ on $L(2c, 1)\cong \mathcal{N}_{c/d}$ for $c>1$ and $\gcd(c,d)=1$, is $r_*\chi$, where $r:ST\S^2\to L(2c,1)$ is the projection.
\end{customthm}
\noindent Note that, for $c>2$, the manifold $L(2c,1)$ is not presented as the unit tangent bundle of a manifold.
Moreover, the study of $\S^2\times \S^1$ allows us to prove an analogue of Theorem \ref{Thm: B} for the class of compact spacetimes $\{(Z\times \S^1,g_Z-\frac{d^2}{c^2}dt^2)\}_{c,d\in \mathbb{N}^+},$ with $Z$ a Zoll manifold (Proposition \ref{prop:Zoll-S1}).

Secondly, we bring methods of Engel geometry to deal with the spaces of null geodesics of three-dimensional spacetimes. An Engel structure on a four-dimensional manifold $Q$ is a rank-two distribution $\mathcal{D}$ that generates a rank-three distribution $\mathcal{E}:=[\mathcal{D}, \mathcal{D}]$ satisfying $[\mathcal{E}, \mathcal{E}] = TQ$, where we are referring to the bracket of sections. The distribution $\mathcal{D}$ defines a unique line distribution $\mathcal{W}$, known as kernel, by the property $[\mathcal{W},\mathcal{E}]\subseteq \mathcal{E}$, which completes a flag $\mathcal{W}\subset\mathcal{D}\subset\mathcal{E}\subset TQ$. There exists a canonical (Cartan) prolongation from a contact three-manifold $(N,\xi)$ to obtain an Engel distribution on $\S(\xi)$. Similarly, given a three-dimensional Lorentzian manifold $(M,g)$, one can canonically define an Engel structure on the projectivization of the bundle of null vectors $\mathcal{P}C$ of $M$. In Proposition \ref{prop:kernel}, we compute the kernel of the Lorentz prolongation of an arbitrary Lorentzian three-manifold. We have later learnt that this computation is also made,  with different techniques, in the preprint \cite[Thm. 1.3]{Mitsumatsu}, but with no mention to the contact structure (which we discuss in Theorem \ref{Thm: C} below).

Next, by considering the Cartan deprolongation of the Lorentz prolongation of a spacetime, we show,
\begin{customthm}{\ref{Thm: C}}\label{Thm: introC}
Let $M$ be a three-dimensional spacetime. Then,
\[
\mathcal{N}  \cong \mathcal{P}C/\mathcal{W}.
\] In addition, if $\mathcal{N}$ is a manifold and $p:\mathcal{P}C\to\mathcal{P}C/\mathcal{W} \cong \mathcal{N}$ is a submersion, the canonical contact structure on $\mathcal{N}$ is
\[
\mathcal{H} \cong p_*\mathcal{E}.
\]
\end{customthm}

Finally, we make use this theorem to make a first step towards the characterization of the three-dimensional contact manifolds that are the space of null geodesics of a spacetime by the existence of a certain Engel manifold together with a foliation, and retrieve the spacetime.

\begin{customthm}{\ref{thm: characterisationEngelFoliation}}\label{Thm: IntroCharacterisation}
A three-dimensional contact manifold $(N, \xi)$ is contactomorphic to the space of null geodesics of a spacetime if and only if there exists an Engel manifold $(Q,\mathcal{D})$ with Engel flag $\mathcal{W}\subset\mathcal{D}\subset\mathcal{E}\subset TQ$ such that
\begin{equation}
N\cong Q/\mathcal{W} \hspace{1cm}\text{and}\hspace{1cm}\xi\cong p_*\mathcal{E},
\end{equation}
for $p:Q\to Q/\mathcal{W}$ the projection, and $Q$ admits an oriented foliation by circles $\mathcal{F}$ such that
\begin{enumerate}
    \item for all $S\in\mathcal{F}$ and $x\in S$, we have $T_xS\oplus \mathcal{W}_x = \mathcal{D}_x$,
    \item the space of leaves $M := Q/\mathcal{F}$ is a manifold and the projection $q:~Q\to M$ is a submersion,
    \item for every $S\in\mathcal{F}$, the image $q_*\mathcal{D}|_S$ is a cone in the vector space $T_{q(S)}M$ and the map $x\in S\mapsto q_*\mathcal{D}_x$ is injective.
\end{enumerate}

In addition, if the above conditions are satisfied, $(N, \xi)$ is contactomorphic to the space of null geodesics of $(Q/\mathcal{F}, g)$, where $g$ is a metric on $Q/\mathcal{F}$ with bundle of cones $q_*\mathcal{D}$.
\end{customthm}

The main objects at play and our results are visually presented as follows.
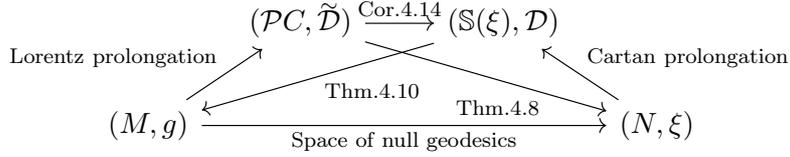
\begin{equation*}
\begin{tikzcd}[column sep=small]
& (\mathcal{P}C,\Dtilde )\arrow[drrr, swap,"\text{Thm.} \ref{Thm: C}", near end] \arrow[rr,rightarrow, "\text{Cor.}\ref{cor: LocalDiffeo}"]
&  & (\S(\xi),\mathcal{D})\arrow[dlll,"\text{Thm.} \ref{thm: characterisationEngelFoliation}"]& \\
(M,g) \arrow[ur, "\text{Lorentz prolongation}"] \arrow[rrrr, swap,"\text{Space of null geodesics}"] & & & & (N,\xi) \arrow[ul,swap, ,"\text{Cartan prolongation}"]
\end{tikzcd}
\end{equation*}
Corollary \ref{cor: LocalDiffeo} deals with the relation between the Cartan and Lorentz prolongations, and opens the interesting question of the relation between the several Engel manifolds having the same Cartan deprolongation.


\medskip

\textbf{Acknowledgments: } We are indebted to Francisco Presas for proposing the ideas that led to this work and for his invaluable comments. We also thank Miguel Sánchez Caja for his interest in this project.

\section{Definitions and basic properties}

\subsection{The space of null geodesics of a spacetime}
We work throughout the paper in the category of smooth manifolds. For a Lorentzian manifold $(M, g)$, a nonzero vector $v\in TM$ is said to be timelike, spacelike or null if $g(v,v)$ is, respectively, negative, positive or zero. 
A smooth curve $\gamma:I\to M$ is timelike, spacelike or null if its velocity vector $\dot{\gamma}$ is so everywhere. Likewise, we talk about timelike, spacelike and null submanifolds or vector fields. 

The set of null vectors on a Lorentzian manifold $M$ has the structure of a smooth bundle $\pi:C\to M$, whose fibres  consist of two hemicones. We denote such a choice by $C^+$. A differentiable choice of one of such hemicones, when possible, makes $M$ time-oriented. Time-orientability is equivalent to the existence of a global timelike vector field $X\in\mathfrak{X}(M)$.

\begin{Definition}
A spacetime is a time-oriented connected Lorentzian manifold of dimension $\geq 3$.
\end{Definition}

\begin{Definition}
\label{Def: SpaceofNullGeod}
The space of null geodesics $\mathcal{N}$ of a spacetime $(M,g)$ is
\[
\mathcal{N} := \{\gamma(I)\ |\ \gamma:I\to M \text{ is a maximal geodesic with }\dot{\gamma}\subset C^+\}.
\]
\end{Definition}

The space $\mathcal{N}$ can be constructed as the leaf space for a distribution on $C^+$. Recall that the geodesic spray $X_g\in\mathfrak{X}(TM)$ is the vector field on $TM$ whose integral lines are $\dot{\gamma}(t)\in T_{\gamma(t)}M$ for $\gamma: I\rightarrow M$ a geodesic, whereas  the Euler vector field $\Delta \in\mathfrak{X}(TM)$ is defined as
\[
\Delta(v) = T_0c(\partial_s),
\]
with $v\in T_xM$ and $c:\R\to T_xM$ given by $c(s) =e^s v$, and whose differential at~$0$ we denote by $T_0c$. Note that $c$ is an integral line of $\Delta$. The geodesic spray and the Euler field are tangent to the bundle $C^+$ and define an integrable distribution $\lspan X_g,\Delta\rspan$ (see \cite{Bautista2008}). Note that, by quotienting $C^+$ by the Euler field, we obtain the projectivization of the bundle $C^+$, which is relevant as we only care about unparametrized null geodesics. Then, by quotienting by the geodesic spray, we identify directions in different projectivized cones for which there exists a geodesic in $M$ going through both of them. Thus, we have 
\begin{equation}\label{eq:Ncong C+ dist}
\mathcal{N} \cong C^+/\lspan X_g,\Delta\rspan.    
\end{equation}

From now on, we will consider the case in which $\mathcal{N}$ is a 
manifold.

\begin{Remark}
A sufficient condition for $\mathcal{N}$ to be a 
manifold is found in \cite{Bautista2008} (namely, when the spacetime $(M, g)$ is strongly causal and null-pseudo-convex).
\end{Remark}

\subsection{The canonical contact structure} Recall that a contact structure on a $(2n+1)-$manifold $N$ is a codimension-one distribution $\xi\subset TN$ which is given, at least locally, as the kernel of a one-form $\alpha$ satisfying that the top form \( \alpha\wedge(d\alpha)^n \)
vanishes nowhere. We give two examples that will be relevant later.

\begin{Example}\label{ex:contact-unit-cotangent}
For any manifold $M$, its unit cotangent bundle $\tilde{\pi}:ST^*M\to M$ has a canonical contact structure (see, for instance, \cite[Lem. 1.2.3]{Geiges2013}). A point $\omega\in ST^*M$ may be regarded as a linear form $\tilde{\omega}\in T^*_{\tilde{\pi}(\omega)}M$ up to positive rescaling, which is  determined by the hyperplane $l_\omega = \ker\tilde{\omega}\subset T_{\tilde{\pi}(\omega)}M$.  The canonical contact distribution on $ST^*M$ is 
	\[
	\xi_\omega := \left(T_\omega \tilde{\pi}\right)^{-1}\left(l_\omega\right).
	\] 
\end{Example}

\begin{Example}
The unit tangent bundle $\pi:STM\to M$ of a Riemannian manifold $(M, g)$, which we define as $STM := \{u\in TM\ |\ g(u,u)=1\}$, has a canonical contact structure coming from $\xi$ for $ST^*M$ as in Example \ref{ex:contact-unit-cotangent}. Regard $g$ as a map $TM\to T^*M$ and consider the contact structure 
\[
\chi := \left(g^{-1}\right)_* \xi.
	\]
	Namely, for $u\in STM$ we have
	\begin{align*}
		\chi_u =& \left(\left(g^{-1}\right)_*\xi\right)_u =  \Big(T_{g(u)}\left(g^{-1}\right)\circ\left(T_{g(u)}\tilde{\pi}\right)^{-1}\Big)\left(l_{g(u)}\right)   \\ = & \left(T_{u}(g\circ\tilde{\pi}
		)\right)^{-1}\big(\ker g(u)\big)=  (T_u\pi)^{-1}(\lspan u\rspan^\perp),
	\end{align*}
where $\lspan u\rspan^\perp$ denotes the orthogonal subspace to $u$ in $T_{\pi(u)}M$ with respect to~$g$. 
\end{Example}

For the smooth manifold $\mathcal{N}$, there is a canonical contact structure defined in terms of the so-called skies of the spacetime.
\begin{Definition}
Let $(M, g)$ be a spacetime and $x\in M$. The sky of $x$ is
\[
\mathfrak{S}_x = \{\gamma\in\mathcal{N}\ |\ x\in\gamma\subset M\}.
\]
\end{Definition}

Note that, for any $x\in M$, the sky $\mathfrak{S}_x$ is in correspondence with the projectivization of the cone $C_x$. Hence,  if $m = \dim M$, we have \(
\mathfrak{S}_x \cong \S^{m-2}.
\)

\begin{Definition}
The canonical contact structure on the manifold of null geodesics $\mathcal{N}$ is the codimension-one distribution $\mathcal{H}$ defined as follows. For $\gamma\in \mathcal{H}$, let $x,y\in\gamma$ such that they cannot be joined by a one-parameter family of geodesics. Then,
\[
\mathcal{H}_\gamma = T_\gamma\mathfrak{S}_x\oplus T_\gamma\mathfrak{S}_y.
\]
\end{Definition}

Note that the existence of two such points $x$ and $y$ follows from the fact that the exponential map provides a local diffeomorphism for any point of~$M$.

\begin{Proposition}\cite[Sec. 2.4]{Bautista2008}
The distribution $\mathcal{H}$ on the manifold $\mathcal{N}$ is well defined and is indeed a contact structure.
\end{Proposition}

\section{The contact manifold of null geodesics of $(\S^2\times\S^1, g_c)$}\label{sec:contact-null-S2S1}

\subsection{The space of null geodesics}\label{sec:null}

Let $(\S^2, g_\circ)$ be the two-sphere with the round metric. Consider $M = \S^2\times\S^1$ and let $t$ be the angle coordinate on $\S^1$. For $x\in \S^2$, we will refer to the point $(x,t)\in M$. Define,  for $c,d\in \mathbb{N}^+$ coprime, the Lorentzian metric 
\[
g _{c/d}= g_\circ -\frac{d^2}{c^2}dt^2.
\]

The pair $(M,g_{c/d})$ is a Lorentzian manifold in which $\S^2\times\{t\}$ is a spacelike surface for any $t\in\S^1$, whereas $\{x\}\times\S^1$ is a  timelike submanifold for any $x\in \S^2$. The vector field $(0,\partial_t)\in T\S^2\oplus T\S^1$ as a choice of future turns $M$ into a spacetime.

\begin{Lemma}\label{lemma:unit-speed}
The space of null geodesics $\mathcal{N}_{c/d}$ of $(M,g_{c/d})$ is given by
\[
\mathcal{N}_{c/d} \cong\left\{\big(\mu(s), \frac{c}{d}s\big)\ |\ \mu \text{ is a unit-speed great circle in }\S^2\right\}.
\]
\end{Lemma}
\begin{proof}
In a product chart using the coordinate $t$ for $\S^1$, the Christoffel symbol $\Gamma_{ij}^k$ of the metric $g_{c/d}$ vanishes whenever $i$, $j$ or $k$ equals 3, and all the others are the Christoffel symbols of $g_\circ$ in the chart of $\S^2$. Hence, the geodesic equation for a curve $\gamma:I\to  M $ defined by $\gamma(s) = \big(\mu(s), t(s)\big)$ is given by  $\ddot{t} = 0$, that is, $t(s) = a+bs$ for some $a, b\in\R$, and the geodesic equation for $\mu$ in $\S^2$. Let $u(s)\in T\S^2$ be the vector tangent to the curve $\mu(s)$. Since $\S^2\times\{t\}$ is a spacelike surface for all $t\in\S^1$, we can suppose, by reparametrizing $\gamma$, that \(
g_{c/d}\big((u,0), (u,0)\big) = g_\circ(u,u) = 1. 
\) Then, 
\[
g_{c/d}(\dot{\gamma}, \dot{\gamma}) = g_\circ(u,u)-\frac{b^2d^2}{c^2} \langle \partial_t, \partial_t\rangle = 1-\frac{b^2d^2}{c^2},
\]
so $\gamma$ is a future-pointing null geodesic if and only if $b = \frac{c}{d}$, as $-\frac{c}{d}$ would give a past-pointing geodesic. By uniqueness of the geodesics in a pseudo-Riemannian manifold, all the null geodesics of $(M, g_{c/d})$ modulus reparametrization are of the~form
\[
\gamma(s) = \big(\mu(s), a+\frac{c}{d}s\big),
\]
where $\mu$ is a unit-speed great circle in $\S^2$. Since $\gamma$ intersects $\S^2\times\{0\}$ at least at one point, we can suppose $a = 0$. \end{proof}

\begin{Remark}
    We use the notation $\mathcal{N}_{c}:=\mathcal{N}_{c/1}$.
\end{Remark}

We start with the case $c = 1$. The speeds at which a geodesic travels the time direction $\S^1$ and the great circle in $\S^2$ are $1$ to $d$.  For $\gamma\in\mathcal{N}_1$, there is a unique $x\in\S^2$ such that $(x,0)\in\gamma$. Indeed, $\gamma$ intersects $\S^2\times \{0\}$ for $s\in 2\pi d\Z$, and $\mu(2\uppi d\Z)$ is the unique point $x$ on $\S^2$. So $\gamma$ is completely determined by $x\in\S^2$ and the tangent vector of the projection $\pi_{\S^2}(\gamma)$ at $x$, which is unitary by Lemma \ref{lemma:unit-speed}. In addition, any $u\in ST\S^2$ defines a unique null geodesic $\gamma\subset M$, which is the lift of the great circle $\mu\subset \S^2$ defined by $u$, meaning that $\pi_{\S^2}(\gamma) = \mu$.
\begin{Proposition}
\label{Prop: STS2}
Thus, we have
\[
\mathcal{N}_{1/d}\cong \mathcal{N}_1 \cong ST\S^2.
\]
\end{Proposition}

\begin{figure}[H]
        \centering
        \resizebox{!}{5cm}{%
            \trimbox{0cm 0.5cm 0cm 0cm}{
           \input{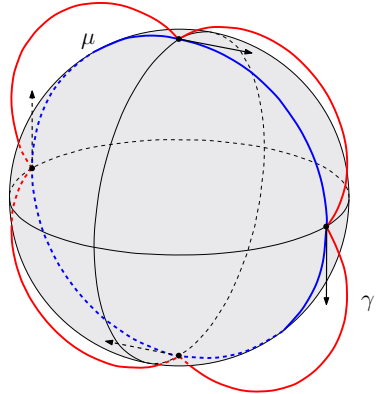}
        }
        }
        \caption{Null geodesic in $\S^2\times\S^1$ with $c = 4$, $d =1$. The grey surface represents $\S^2\times\{0\}$ and the radial coordinate is the $\S^1$ direction. The four elements of $ST\S^2$ represent the same red null geodesic $\gamma$.}
        \label{Fig: cmajor1}
    \end{figure}

Let us now consider $\mathcal{N}_{c/d}$ with $c>1$. For a geodesic, the ratio between the turns around the time direction $\S^1$ and the turns around a great circle in $\S^2$ is $c$ to $d$. Since geodesics are travelled at constant speed, every $\gamma\in\mathcal{N}_{c/d}$ intersects $\S^2\times\{0\}$ at $c$ points, namely those in $\mu(2\uppi\frac{d}{c}\Z) = \mu(\frac{2\uppi}{ c}\Z)$, where the equality follows from the fact that $\gcd(c,d)=1$. These points are equidistantly spread over the great circle $\mu := \pi_{\S^2}(\gamma)$, see Figure \ref{Fig: cmajor1}.   Conversely, any $u\in ST\S^2$ defines a unique null geodesic, which is the lift of the great circle defined by $u$ to $M$.

In order to get a proper description of $\mathcal{N}$, one ought to identify the different elements of $ST\S^2$ defining the same null geodesic. If $\gamma$ is a null geodesic of $M$ that intersects $\S^2\times\{0\}$ at $(x,0)$ with tangent vector $(u, \partial_t)\in ST\S^2\oplus T\S^1$, then $\gamma$ intersects $\S^2\times\{0\}$ at $(x_j, 0): = \left(
\mu\left(\frac{2\uppi dj}{c}\right), 0\right)$ for $j = 0,\ldots,c-1$, where $\mu\subset\S^2$ is the great circle defined by $u$, and with velocity $(u_j,\partial_t):=\left(\dot{\mu}\left(\frac{2\uppi d j}{c}\right), \partial_t\right)$. Note that $(x_j,u_j)$ can be obtained by a rotation of $(x,u)$ of $\frac{2\uppi dj}{c}$ radians about the axis $x\times u$, 
\begin{equation*}
\begin{pmatrix}x_j\\ u_j
\end{pmatrix}
= \begin{pmatrix}
	\cos\frac{2\uppi dj}{c} & \sin\frac{2\uppi dj}{c}\vspace{0.08cm}\\
	-\sin\frac{2\uppi dj}{c}&
	\cos\frac{2\uppi dj}{c}
\end{pmatrix}\begin{pmatrix}
x\\u
\end{pmatrix}.
\end{equation*}

Hence, we have shown, 

\begin{Proposition}
\label{Prop: Zc}
Consider the $\Z_c$-action on $ST\S^2$ generated by 
\begin{equation*}
	\label{ActionZc}
\begin{pmatrix}y\\ v\end{pmatrix}\mapsto \begin{pmatrix}
	\cos\frac{2\uppi d}{c} & \sin\frac{2\uppi d}{c}\vspace{0.08cm}\\
	-\sin\frac{2\uppi d}{c}&
	\cos\frac{2\uppi d}{c}
\end{pmatrix}\begin{pmatrix}
	y\\v
\end{pmatrix}.
\end{equation*}

Then,
\[
\mathcal{N}_{c/d}\cong ST\S^2/\Z_c,
\]
where $ST\S^2/\Z_c$ denotes the orbit space of $ST\S^2$ under the action of $\Z_c$.
\end{Proposition}


The next step is to obtain an explicit description of the spaces $\mathcal{N}_{c/d}$ that allows us to compute their canonical contact structure. We will prove the following result.

\begin{Theorem}
\label{Thm: A}
For any $c,d\geq 1$ with $\gcd(c,d)=1$,
\[
\mathcal{N}_{c/d}\cong L(2c,1).
\]
\end{Theorem}

\subsection{A quaternionic approach to $ST\S^2$ and the Hopf fibration}\label{sec:quaternionic}
In order to prove Theorem \ref{Thm: A}, we develop a quaternionic approach to $ST\S^2$ and the lens spaces $L(2c,1)$. Let $\H$ denote the division algebra of quaternions and $\V$ the three-dimensional vector space of pure imaginary quaternions. The canonical identification of $\H$ with $\R^4$ defined by $\alpha+ai+bj+ck\mapsto(\alpha,a,b,c)$ gives an identification $\V\cong \R^3$, which provides $\V$ with a cross product induced by that of $\R^3$: for $u$, $v\in\V$,
\[
u\times v = \frac{uv-vu}{2}.
\]

Let $^*:\H\to\H$ be the conjugation on $\H$, which is an antiautomorphism allowing us to define a norm  $|q|^2 = qq^*$ for $q\in\H$. The restriction of such norm on $\V$ induces, via the polarization identity, an inner product on $\V$ defined, for $u$, $v\in\V$, by
\[
\langle u,v\rangle = -\frac{uv+vu}{2},
\]
which coincides with the euclidean inner product in $\R^3$. We can also identify $\S^3\cong S\H := \{q\in\H\ |\ |q| = 1\}$ and $\S^2\cong S\V : = \{u\in\V\ |\ \langle u,u\rangle = -u^2 = 1\}$. Finally, one has $ST\S^2\cong ST(S\V):=\{(u,v)\in S\V\times S\V\ |\ \langle u,v\rangle = 0\}$.

\begin{Lemma}
\label{Prop: Hopf}
    For any $w\in S\V$, there exists a Hopf-like fibration map
	\[
	\begin{array}{cccc}
		\tau_w:&S\H&\to&S\V\\
		& q&\mapsto&q^{-1}wq,
	\end{array}
	\]
	which provides $S\H$ with the structure of an $\S^1$-bundle over $S\V$. The fibre over $p\in S\V$ is given by $\{e^{w\theta}q\ |\ \theta\in\R\}$, for any $q\in\tau_w^{-1}(p)$.
\end{Lemma}
\begin{proof}
Let $w\in S\V$. Let us show first that $\tau_w$ maps onto $S\V$. Indeed, for $q\in S\H$, 
\[
	q^{-1}wq+(q^{-1}wq)^* = q^{-1}wq-q^{-1}wq = 0,
\]
\[
	\langle q^{-1}wq, q^{-1}wq\rangle = -(q^{-1}wq)(q^{-1}wq) = -q^{-1}w^2q = q^{-1}q = 1.
\]

We show next that the map is surjective. Let $p\in S\V$ and assume it is not collinear with $w$. Let $\theta = \arccos\langle w, p\rangle\in (0,\uppi)$ and $\eta= \frac{p\times w}{|p\times w|}\in S\V$. Define $q := e^{\eta \frac{\theta}{2}} = \cos\frac{\theta}{2}+\eta \sin\frac{\theta}{2}$. We have $ q^{-1}wq=(\cos \theta ) w + (\sin \theta) (w\times \eta)$, 
a rotation of $w$ of angle $-\theta$ around the axis $\eta$. 
Hence, $q^{-1}wq = p$, as needed. If $p$ and $w$ are collinear, take $q = 1$ or $q = e^{u\frac{\uppi}{2}}$, with $u\in S\V$ perpendicular to $w$, depending on whether $p = w$ or $p= -w$.

Finally, let $q\in \tau_w^{-1}(p)$ be arbitrary. Then, if $\theta\in \R$, 
\[
\tau_w(e^{w\theta}q) = q^{-1}e^{-w\theta}we^{w\theta}q = q^{-1}wq = \tau_w(q).
\]
In addition, if $\tau_w(q_1) = \tau_w(q_2)$, then $w = (q_2q_1^{-1})^{-1}w(q_2q_1^{-1})$, which implies that $q_2q_1^{-1} = e^{-w\theta}$ for some $\theta\in\R$.
\end{proof}

\begin{Proposition}\label{prop:Phi:S3-STS2}
The map
\[
\begin{array}{cccc}
     \Phi:& S\H & \to & ST(S\V) \\
     & q&\mapsto & (q^{-1}kq, q^{-1}jq)
\end{array}
\]
provides a surjective local diffeomorphism in such a way that the preimage of a point in $ST(S\V)$ consists of exactly two antipodal points in $S\H$.
\end{Proposition}
\begin{proof}
	Surjectivity follows using the same ideas as in the proof of Lemma \ref{Prop: Hopf}. Also, for any $q\in S\H$, we have $\Phi(q) = \Phi(-q)$. If $q_1, q_2\in S\H$ are such that their images under $\Phi$ coincide, then $q_2 q_1^{-1} = e^{k\theta_k} = e^{j\theta_j}$ for some $\theta_j, \theta_k\in[0,2\uppi)$, which can only happen if $\theta_k = \theta_j = 0$ or $\theta_k = \theta_j = \uppi$.
\end{proof}

\begin{Remark}
Actually, in Proposition \ref{prop:Phi:S3-STS2} it is possible to replace $j$ and $k$ by any $u$, $v\in S\V$ such that $\langle u, v\rangle = 0$ and the result remains true.
\end{Remark}

\begin{Corollary}
\label{Cor: SHZ2}
	Let $\Z_2$ act on $S\H$ via the antipodal map. There exist diffeomorphisms $S\H/\Z_2\cong ST(S\V)$, and hence $\S^3/\Z_2\cong ST\S^2$, which we also denote by~$\Phi$. 
\end{Corollary}

We establish now the connection to lens spaces.

\begin{Definition}
\label{Def: Lens}
Consider the 3-sphere $\S^3\subset\C^2$ and $p\in\Z^+$. Define the $\Z_p-$action on $\S^3$ generated by
\[
(z_0,z_1)\mapsto \left(e^{\frac{2\uppi i }{p}}z_0, e^{\frac{2\uppi i }{p}}z_1\right).
\] The lens space $L(p,1)$ is the smooth manifold $L(p,1) := \S^3/\Z_p$.
\end{Definition}

Identifying $\R^4\cong\C^2$ with $\H$ via $(z_0,z_1)\mapsto z_0+z_1j$, the $\Z_p$ action in Definition~\ref{Def: Lens} becomes $q\mapsto e^{\frac{2\uppi i}{p}}q$. Then, $L(p,1)\cong S\H/\Z_p$. 

Since the $\Z_2$-action on $S\H$ that defines the lens space $L(2,1)$ is precisely the given by the antipodal map, we have shown,

\begin{Proposition}
\label{Prop: LensSTS2}
	We have $ST(S\V)\cong S\H/\Z_2\cong L(2,1)$, so $ST\S^2 \cong L(2,1)$.
\end{Proposition} 

This quaternionic approach allows us to formalize and prove the following result, which will give us, together with Propositions \ref{Prop: STS2} and \ref{Prop: Zc}, the proof of Theorem ~\ref{Thm: A}.

\begin{Proposition}
\label{Prop: LensZc}
	Let $\Phi:S\H/\Z_2\to ST(S\V)$ and $c\geq 2$. Then, the $\Z_{2c}$-action on $S\H$ that generates the lens space $L(2c,1)$ descends to a $\Z_c$-action on $S\H/\Z_2$ that, via $\Phi$, induces the $\Z_c$-action on $ST(S\V)$ generated by 
	\[\begin{pmatrix}u\\v\end{pmatrix} \mapsto \begin{pmatrix}
		\cos\frac{2\uppi}{c} &  	\sin\frac{2\uppi}{c}\vspace{0.08cm}\\
		-\sin\frac{2\uppi}{c}&  	\cos\frac{2\uppi }{c}
	\end{pmatrix}\begin{pmatrix}
		u\\v
	\end{pmatrix}.\]
	Hence, $\Phi$ induces a diffeomorphism between $S T(S\V)/\mathbb{Z}_c$ and $L(2c,1)$.
\end{Proposition}
\begin{proof}
	
	Let $q\in S\H$. The $\Z_{2c}$-action on $S\H$ is generated by $q\mapsto e^{\frac{\uppi i}{c}}q$. Then, 
	\[
	\Phi(e^{\frac{\uppi i}{c}}q) = \Big(q^{-1}e^{-\frac{\uppi i}{c}}ke^{\frac{\uppi i}{c}}q, q^{-1}e^{-\frac{\uppi i}{c}} j e^{\frac{\uppi i}{c}}q\Big) = \Big(q^{-1}ke^{\frac{2 \uppi i}{c}}q, q^{-1}je^{\frac{2 \uppi i}{c}}q\Big).
	\]
	
	Now, consider $\Phi(q) = (q^{-1}kq, q^{-1}jq)$. The result follows from 
	\begin{align*}
		\cos\frac{2\uppi}{c}q^{-1}kq+\sin\frac{2\uppi}{c}q^{-1}jq & {} = q^{-1}k\Big(\cos\frac{2\uppi}{c}+i\sin\frac{2\uppi}{c}\Big)q = q^{-1}ke^{\frac{2\uppi i }{c}}q,\\
		-\sin\frac{2\uppi}{c}q^{-1}kq+\cos\frac{2\uppi}{c}q^{-1}jq & {} = q^{-1}j\Big(\cos\frac{2\uppi}{c} + i\sin\frac{2\uppi}{c}\Big)q = q^{-1}je^{\frac{2\uppi i}{c}}q.
	\end{align*}\end{proof}

Proposition \ref{Prop: LensZc} is enough to prove Theorem \ref{Thm: A} for $d = 1$. For $d>1$, we need the following observation.

\begin{Lemma}\label{Lem: GroupAuto}
  Let $c,d\in \Z^+$ coprime. The group automorphism of $\Z_{c}$ defined by $m\mapsto d\cdot m$ sends the $\Z_c-$action on $ST\S^2$ in Proposition \ref{Prop: LensZc} to the $\Z_c-$action given by Equation \eqref{ActionZc}.
 Hence, the two orbit spaces are diffeomorphic.
\end{Lemma}

\begin{proof}[Proof of Theorem \ref{Thm: A}]
The case $d,c = 1$ follows from Proposition \ref{Prop: STS2} and Proposition \ref{Prop: LensSTS2}. The rest of the cases with $d = 1$ follow from Proposition \ref{Prop: Zc} and Proposition \ref{Prop: LensZc}. For $d>1$, we make use of Lemma \ref{Lem: GroupAuto}.
\end{proof}

\subsection{The canonical contact structures}\label{sec:canonical-contact}
We next compute the canonical contact structures on the spaces $L(2c,1)$, $c\geq 1$, which were seen as spaces of null geodesics. Note that, since the construction of the canonical structure is completely local, we can assume, without loss of generality, that $d = 1$. We do this throughout.

\begin{Theorem}
\label{Thm: B}
The canonical contact structure $\mathcal{H}$ on $ST\S^2\cong\mathcal{N}_{1/d}$ is the canonical contact structure $\chi$ on $ST\S^2$. In general, the canonical contact structure $\mathcal{H}$ on $L(2c, 1)\cong \mathcal{N}_{c/d}$ for $c>1$ and $\gcd(c,d)=1$, is $r_*\chi$, where $r:ST\S^2\to L(2c,1)$ is the projection.
\end{Theorem}

 We break the proof of Theorem \ref{Thm: B} into two lemmas, proving first the case $c = 1$ and then the cases with $c>1$.

\begin{Lemma}\label{lem:can-contact-on-STS2}
The canonical contact structure $\mathcal{H}$ on $ST\S^2\cong\mathcal{N}_1$ is the canonical contact structure $\chi$ on $ST\S^2$.
\end{Lemma}
\begin{proof}
Let $\gamma\in\mathcal{N}_1\cong ST\S^2$. Recall that $\gamma$ is the lift of the great circle $\mu:\R\to \S^2$ defined by the pair $(x,u)\in ST\S^2$ representing $\gamma$. We will show that $\mathcal{H}_\gamma = \chi_{(x,u)}$.

Take $(x,0)\in\gamma$ and $\gamma(\tau)\neq (x,0)$, with $0<\tau<\uppi$. Note that all geodesics in the sky $\mathfrak{S}_{(x,0)}$ intersect $\S^2\times\{0\}$ at $(x,0)$, so its projection is $\pi(\mathfrak{S}_{(x,0)}) = \{x\}$. Hence, $T_{x}\pi(T_\gamma\mathfrak{S}_{(x,0)}) = \{ 0\} \subseteq \lspan u\rspan^\perp$ and $T_\gamma\mathfrak{S}_{(x,0)}\subset \chi_{(x,u)}$.

Consider now $\gamma(\tau)\neq (x,0)$. Since $c = 1$, we know that $\gamma(s) = \big(\mu(s), s\big)$ and hence $ \pi_{\S^2}\big(\gamma(\tau)\big) = \mu(\tau)=:y$. Let $v\in ST_y\S^2$ such that $\langle\dot{\mu}(\tau), v\rangle = 0$. Since all geodesics in $\S^2\times\S^1$ are travelled at the same speed, the projection of the sky of $\gamma(\tau)$ is parametrized by $\pi(\mathfrak{S}_{\gamma(\tau)})(s) = y\cos\tau+(\dot{\mu}(\tau)\cos s+v\sin s)\sin\tau$, which is a circle on $(\S^2,g_\circ)$ of radius $\tau$ and centre $y$, see Figure \ref{Fig: circleonS2}, and $\pi(\gamma) = \pi(\mathfrak{S}_{\gamma(\tau)})(0)$. Hence, 
\begin{align*}
T_x\pi(T_\gamma\mathfrak{S}_{\gamma(\tau)}) &= \lspan\frac{d}{ds}\Big|_{s = 0}\big(y\cos\tau+(\dot{\mu}(\tau)\cos s+v\sin s)\sin\tau\big)\rspan \\ &= \lspan v\sin\tau\rspan\subset\lspan u\rspan^\perp,
\end{align*}
which implies $T_\gamma\mathfrak{S}_{\gamma(\tau)}\subset \chi_{(x,u)}$.

  \begin{figure}[H]
        \centering
        \resizebox{!}{5cm}{%
            \trimbox{0cm 8cm 0cm 9cm}{
               \input{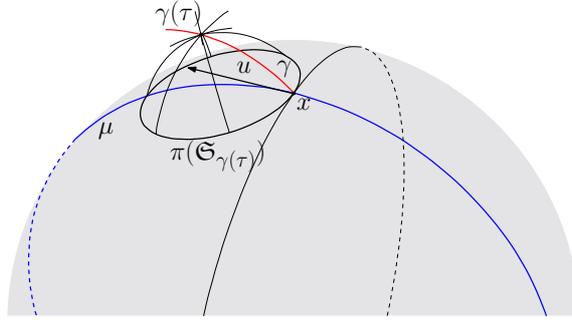}
        }
        }
        \caption{The projection onto $\S^2$ of the sky of $\gamma(\tau)$ is a circle tangent to $\lspan u\rspan^\perp$.}
        \label{Fig: circleonS2}
    \end{figure}

Since the distributions $\mathcal{H}$ and $\chi$ on $ST\S^2$ both have rank 2, the Lemma is proved.
\end{proof}

Let $r:ST\S^2\cong L(2,1)\to L(2c,1)$ be the canonical projection for $c>1$. We use the notation $[u]:=r(u)\in L(2c, 1)$, which is the class of $u\in ST\S^2$ under the action of $\Z_c$. Let $[u]\in L(2c,1)$ and $U$ a neighbourhood of $u$ for which $r|_{U}:U\to r(U)$ is a diffeomorphism. We will show that $$\mathcal{H}_{[u]} = (r|_U)_*\chi_{u}.$$

The following lemma concludes the proof of Theorem \ref{Thm: B}.
\begin{Lemma}\label{lem:can-contact-on-STS2-2}
For $c>1$, the canonical contact structure $\mathcal{H}$ on $L(2c,1)\cong \mathcal{N}_c$, as the space of null geodesics, is $r_*\chi$.
\end{Lemma}
\begin{proof}
Let $x := \pi(u)$. We know that $[u]$ describes the geodesic $\gamma$ in $\S^2\times\S^1$ that intersects $\pi(U)\times\{0\}$ only at $(x,0)$. Take $(x,0)\in\gamma$ and consider its sky $\mathfrak{S}_{(x,0)}$. It is clear that $\mathfrak{S}_{(x,0)} = \{[v]\ |\ v\in ST_{x}\S^2\}$, and thus
\begin{align*}
	\big(\pi\circ r|_{U}^{-1}\big)_* \big(T_{[u]}\mathfrak{S}_{(x,0)}\big) = & T_{x}\left(\left\{\pi(v)\ |\ v\in ST_{x}\S^2\right\}\right) \\ = & T_{x}(\{x\}) = \{0\}\subset\lspan u\rspan^\perp,
\end{align*}
from which we deduce that $(r|_{U}^{-1})_* \big(T_{[u]}\mathfrak{S}_{(x,0)}\big)\subset \chi_{u}$ and $T_{[u]}\mathfrak{S}_{(x,0)}\subset (r|_U)_*\chi_{u}.$

Take now $\gamma(\tau)\neq(x,0)$ but close enough so that $\mathfrak{S}_{\gamma(\tau)}\subset U$. As discussed previously, $\pi\circ  r|_U^{-1}(\mathfrak{S}_{\gamma(\tau)})$ describes a circle $\phi(s)$ in $\pi(U)$ whose tangent vector at $x$ is orthogonal to $u$. Hence,
\begin{align*}
	(\pi\circ r|_U^{-1})_* \big(T_{[u]}\mathfrak{S}_{\gamma(\tau)}\big) =& T_{x}\Big(\pi\circ r|_U^{-1}(\mathfrak{S}_{\gamma(\tau)} )\Big) \\=& T_{x}\{\phi(s)\ |\ s\in \mathbb{R}\}\subset\lspan u\rspan^\perp,
\end{align*}
which implies that $(r|_U^{-1})_* \big(T_{[u]}\mathfrak{S}_{\gamma(\tau)}\big)\subset \chi_{u},$ and hence $T_{[u]}\mathfrak{S}_{\gamma(\tau)}\subset (r|_U)_*\chi_{u}.$

Therefore, $\mathcal{H}_\gamma = T_{[u]}\mathfrak{S}_{(x,0)}\oplus T_{[u]}\mathfrak{S}_{\gamma(\tau)} \subset (r|_U)_*\chi_u$ and, since $r|_U$ is a submersion, equality follows. If $v\in ST\S^2$ is such that $r(v) = [u] = r(u)$ and $V$ is a neighbourhood of $v$ such that $r|_V$ is a diffeomorphism, following the same argument, we have $(r|_V)_*\chi_v = \mathcal{H}_\gamma = (r|_U)_*\chi_u$, and $r_*\chi$ is indeed well defined.
\end{proof}

\subsection{Generalization to a class of compact spacetimes}\label{sec:generalization-compact}

We extend now the results of Section  \ref{sec:canonical-contact} to the class of spacetimes $$\{(Z\times \S^1,g_Z-\frac{d^2}{c^2}dt^2)\}_{c,d\in \mathbb{N}^+ \text{ coprime}},$$ where $(Z,g_Z)$ is a Zoll manifold, that is, such that all of its geodesics are closed and of the same minimal period (which we normalize to $2\uppi$). By Hopf-Rinow theorem, every Zoll manifold is compact, since it is clearly bounded.

\begin{Proposition}\label{prop:Zoll-S1}
The space of null geodesics $\mathcal{N}_{c/d}$ of $(Z\times \S^1,g_Z-\frac{d^2}{c^2}dt^2)$ with $c,d$ coprime, is diffeomorphic to the manifold $STZ/\mathbb{Z}_c$, where the generator of $\Z_c$  identifies tangent vectors on the same geodesic after $1/c$ of a turn, with the pushforward of the canonical contact structure on $STZ$.
\end{Proposition}

\begin{proof}
The same argument as in the proof of Proposition \ref{Prop: STS2} implies $\mathcal{N}_{1/d}\cong STZ$. 
For the contact structure, we can assume $d = 1$. We follow the proof of Lemma \ref{lem:can-contact-on-STS2}, which is geometrically more intuitive. Given a geodesic $\gamma$ and $(x,0)\in\gamma\subset Z\times \S^1$, all geodesics in the sky $\mathfrak{S}_{(x,0)}$ intersect $Z\times \{0\}$ at $(x,0)$ and we analogously have $T_\gamma \mathfrak{S}_{(x,0)} \subset \chi_{(x,u)}$. For a different point $\gamma({\tau})$ in the geodesic (close enough to $\gamma(0)$), the sky $\mathfrak{S}_{\gamma({\tau})}$ is a sphere consisting of centre $y:=\pi_Z(\gamma(\tau))$ and radius the distance between $x$ and $y$. Since $u$ points from $x$ towards the direction of $y$, we have in general $T_x\pi_Z(T_\gamma\mathfrak{S}_{\gamma(\tau)})\subset \lspan u \rspan^\perp$. Indeed, let $w\in T_yZ$ such that $x = \exp_y(w)$ and complete it to an orthogonal basis $(w, w_2,\ldots, w_n)$ of $T_yZ$. The sphere $\pi_Z(\mathfrak{S}_{\gamma(\tau)})$ is the image under $\exp_y$ of the sphere in $T_yZ$ of radius $\sqrt{g_Z(w,w)}$. Now, note that $T_w \exp_y(w)\in \lspan u\rspan$, and $ T_x(\pi_Z(\mathfrak{S}_{\gamma(\tau)}))= \lspan T_w \exp_y(w_2), \ldots , T_w \exp_y(w_n)\rspan$. By Gauss's Lemma,
\[
\langle T_w \exp_y (w), T_w \exp_y (w_i)\rangle = \langle w, w_i\rangle = 0
\]
for all $2\leq i\leq n$, and the claim follows.

For the case $c\geq 2$ and $d =1$ we have an analogue of the end of Section \ref{sec:null}. The  generator of $\Z_c$  identifies tangent vectors on the same geodesic after $1/c$ of a turn. The action is then free (as only the identity element would fix a vector). Since $\Z_c$ is a finite group, the action is proper and the space $STZ/\Z_c$ is a manifold. For $d>1$, we have an analogue of Proposition \ref{Prop: LensZc} and the lemma right after ensuring that $\mathcal{N}_{c/d}\cong \mathcal{N}_c$. For the contact structure, the proof of Lemma \ref{lem:can-contact-on-STS2-2} applies by taking $\gamma(\tau)$ close enough to $(x,0)$ so that the action of the group is trivial and replacing the circle $\phi(s)$ by the sphere described above.
\end{proof}

Note that the purpose of Section  \ref{sec:quaternionic} is giving a concrete description of $STZ/\Z_c$ for $c\geq 2$, which we cannot have in the generality of Proposition \ref{prop:Zoll-S1}, in the case of $Z=\S^2$.

The main examples of Zoll manifolds are Zoll surfaces (which are always spheres with Zoll metrics) and compact symmetric spaces of rank one \cite{Besse1978}. 

\section{Engel structures as a tool in retrievability}
\subsection{Engel geometry and prolongations}

We recall here the main definitions 
on Engel manifolds and present the Cartan and Lorentz prolongations of, respectively, a contact and a Lorentzian three-manifold \cite{Casals_2017, del_Pino_2017}.

 A rank-three distribution $\mathcal{E}\subset TQ$ on a four-manifold $Q$ is said to be an even-contact structure if it is everywhere non-integrable, that is, if $[\mathcal{E},\mathcal{E}] = TQ$.

\begin{Definition}
	Let $Q$ be a four-manifold. A rank-two distribution $\mathcal{D}\subset TQ$ on $Q$ is an Engel structure if $\mathcal{E} := [\mathcal{D},\mathcal{D}]$ is an even-contact structure on $Q$.
\end{Definition}

An Engel structure $\mathcal{D}$ on a four-manifold $Q$ defines a unique line field $\mathcal{W}\subset \mathcal{E}$ by the relation $[\mathcal{W}, \mathcal{E}]\subseteq\mathcal{E}$. The line field $\mathcal{W}$ is known as the kernel (or characteristic line field) of the distribution and it can be shown to lie in the two-distribution $\mathcal{D}$, for which it completes a flag $\mathcal{W}\subset\mathcal{D}\subset\mathcal{E}\subset TQ$. 

\begin{Example}[Cartan prolongation]
	\label{Ex: CartanProl}
	Let $(N,\xi)$ be a contact three-manifold and consider the $\S^1-$bundle $\pi_C: \S(\xi)\to N$, where $\S(\xi)_x$ is the quotient of $\xi_x\setminus\{0\}$ by the relation $v\sim\lambda v$ for all $\lambda\in\R^+$. We regard points in $\S(\xi)$ as pairs $(x,R)$ with $x\in N$ and $R$ an oriented line in $\xi_x$. The canonical Engel structure on $\S(\xi)$ is
	\[
	\mathcal{D}_{(x,R)} := (T_{(x,R)}\pi_C)^{-1}(R).
	\]
	Let $\lspan V,Y\rspan = \xi$ be a local frame on an embedded ball $B\subset N$. Then, $B\times~\S^1~\cong~\S(\xi)|_B$ via $(x,t)\mapsto(x,R := \orlspan X:= V\cos t+Y\sin t\orrspan)$, where $\orlspan\phantom{-}\orrspan$ denotes the oriented spanned line. 
	If we let the dot denote derivation with respect to the coordinate on the fibre, then $\mathcal{D}=~\lspan \partial_t,X \rspan$, $\mathcal{E} = \lspan\partial_t,X, \dot{X}\rspan = \lspan\partial_t\rspan\oplus\xi$ and, since $[\partial_t,\ddot{X}] = -\dot{X}\in\mathcal{E}$, we have $\mathcal{W} = \lspan\partial_t\rspan$.
\end{Example}
\begin{Remark} \label{Rk: CartanExample}
Whenever the leaf space $\S(\xi)/\mathcal{W}$ is a manifold, then  $$N \cong \S(\xi)/\lspan\partial_t\rspan = \S(\xi)/\mathcal{W}.$$ Also, if the projection $p:\S(\xi)\to \S(\xi)/\mathcal{W}$ is a submersion, we have $\xi\cong p_*\mathcal{E}$.
\end{Remark}
\begin{Example}[Lorentz Prolongation]
	\label{Ex: Lorentzprol}
	Let $M$ be a Lorentzian three-manifold. The set of null vectors on $M$ induces an $\S^1-$bundle $\pi_L: \mathcal{P}C\to M$, where $\mathcal{P}C$ is fibrewise the projectivization of the cone $C$. A point $(x,l)\in\mathcal{P} C$ consists of a point $x\in M$ and a line $l$ in $C_x$.  Define an Engel structure on $\mathcal{P} C$ at $(x,l)$ by
	\[
	\mathcal{D}_{(x,l)} := (T_{(x,l)}\pi_L)^{-1}(l).
	\]
	Let $(V,Y,T)$ be an orthonormal frame of $TM$, with $V,Y$ spacelike and $T$ timelike, on an open ball $B\subset M$. Then, $B\times\S^1\cong\mathcal{P}C|_B$ via $(x,\theta)\mapsto(x,l: = \lspan X:= V\cos \theta+Y\sin\theta+T\rspan)$, where the vector fields are on $x$. 
	Letting the dot denote derivation with respect to the fibre coordinate, $[\partial_\theta, \dot{X}]\notin\mathcal{E} = \lspan\partial_\theta, X,\dot{X}\rspan$, which implies that $\mathcal{W}$ is always transverse to~$\partial_\theta$.
\end{Example}
Following the ideas of Remark \ref{Rk: CartanExample}, we make the following observation, which will be useful in Section \ref{sec:recovering}.
\begin{Remark}
\label{Rk: LorentzExample} The family of skies $\{\mathcal{P}C_x\}_{x\in M}$ defines a circle foliation of $\mathcal{P}C$ whose leaf space is diffeomorphic to the manifold $M$. In addition, the bundle of null cones of $M$ can be recovered via the pushforward of the Engel distribution $\mathcal{D}$ under the projection map $\pi_L$, that is, $C_x = (\pi_L)_*\big(\mathcal{D}|_{\mathcal{P}C_x}\big)$ for all $x\in M$.
\end{Remark}

From now on, we denote elements of $\mathcal{P}C$ and $\S(\xi)$ by the line or oriented line that they define, dropping the base point of the three-manifold $M$.

\subsection{The space of null geodesics as a deprolongation} We make use of the deprolongation procedure in Remark \ref{Rk: CartanExample} to present the space of null geodesics as a Cartan deprolongation of the Lorentz prolongation.

\begin{equation*}\label{diag:Engelx1}
\begin{tikzcd}[column sep=small]
& (\mathcal{P}C,\mathcal{D} ) \arrow[dr, dashed, "?"]& \\
(M,g) \arrow[ur, "\text{Ex.}\ref{Ex: Lorentzprol}"] \arrow[rr, "\text{Def.}\ref{Def: SpaceofNullGeod}"] &  & (\mathcal{N},\mathcal{H})
\end{tikzcd}
\end{equation*}

We first recall a technical result.

\begin{Proposition}[\cite{Kowalski_2013}]
\label{Prop: DiagonalMetric}
	Any point of a pseudo-Riemannian three-manifold admits a local chart in which the metric is diagonal.
\end{Proposition}

Let $(M,g)$ be a three-dimensional spacetime and consider $x\in M$. Let $\varphi:(x_1,x_2,x_3) \in V\mapsto\varphi\big((x_1,x_2,x_3)\big)\in U$ be local coordinates around $x$ for which $g$ is diagonal. The matrix representation of $g$ in the chart $(U, \varphi^{-1})$ is 
\begin{equation}
	\label{matrixdiag}
	g\big(\varphi(x_1,x_2,x_3)\big) = 
\begin{pmatrix}
	g_{11}(x_1,x_2,x_3) & 0 & 0\\
	0 & g_{22}(x_1,x_2,x_3)  & 0\\
	0 & 0 & g_{33}(x_1,x_2,x_3) 
\end{pmatrix}
\end{equation}
for some smooth functions $g_{11}, g_{22}, g_{33}$ on $V$. Since the metric is non-degenerate at every point, we can assume $g_{11},g_{22}>0$ and $g_{33}<0$. In addition, the coordinate vector fields $u_i :=\varphi_*e_i$ give the eigendirections of the metric at every point. This discussion allows us to define local coordinates on $\mathcal{P}C$ via
\[
\begin{array}{cccc}
	\Psi: & V\times (0,2\uppi) & \to & \Psi\big(V\times (0,2\uppi)\big)\vspace{0.08cm}\\
	&(x_1,x_2,x_3,\theta)  & \mapsto & \lspan\sdfrac{\cos\theta}{\sqrt{g_{11}}}u_1 + \sdfrac{\sin\theta}{\sqrt{g_{22}}}u_2+\sdfrac{1}{\sqrt{-g_{33}}}u_3\rspan\subset T_{\varphi(x_1,x_2,x_3)}M.
\end{array}
\] For the rest of this section, we denote a line in $\mathcal{P}C$ by the vector that spans it, identifying $\mathcal{P}C\cong\{\sdfrac{\cos\theta}{\sqrt{g_{11}}}u_1 + \sdfrac{\sin\theta}{\sqrt{g_{22}}}u_2+\sdfrac{1}{\sqrt{-g_{33}}}u_3\}_{\theta\in\S^1, x\in M}$. Let us denote by $\partial_{x_1}, \partial_{x_2}, \partial_{x_3}, \partial_{\theta}$ the coordinate vector fields defined by $\Psi$. We compute next the kernel of the Lorentz prolongation of any Lorentzian three-manifold. 

\begin{Proposition}\label{prop:kernel}
	\label{PropKernel}
	In the notation above, the kernel $\mathcal{W}$ of the Engel distribution on $\mathcal{P}C$ defined by the Lorentz prolongation is spanned, on $\Psi\big(V\times (0,2\uppi)\big)$, by the vector field
	\[
	Z := \frac{\cos\theta}{\sqrt{g_{11}}}\partial_{x_1} + \frac{\sin\theta}{\sqrt{g_{22}}}\partial_{x_2}+\frac{1}{\sqrt{-g_{33}}}\partial_{x_3}+(F\cos\theta+G\sin\theta+H)\partial_\theta,
	\]
	where we define
		\[
	\begin{cases}
	\begin{aligned}
	F &:=\phantom{-} \sdfrac{1}{2g_{11}}\left(\sdfrac{1}{\sqrt{g_{22}}}\sdfrac{\partial g_{11}}{\partial x_{2}}+\sdfrac{\sin\theta}{\sqrt{-g_{33}}}\sdfrac{\partial g_{11}}{\partial x_3}\right),
		\\
		G &:= -\sdfrac{1}{2g_{22}}\left(\sdfrac{1}{\sqrt{g_{11}}}\sdfrac{\partial g_{22}}{\partial x_{1}}+\sdfrac{\cos\theta}{\sqrt{-g_{33}}}\sdfrac{\partial g_{22}}{\partial x_3}\right),\\
		H & :=\phantom{-}\sdfrac{1}{2g_{33}}\left( \sdfrac{\sin\theta}{\sqrt{g_{11}}}\sdfrac{\partial g_{33}}{\partial x_1}-\sdfrac{\cos\theta}{\sqrt{g_{22}}}\sdfrac{\partial g_{33}}{\partial x_2}\right).
		\end{aligned}
	\end{cases}
	\]
\end{Proposition}
\begin{proof}
	Recall, from Example \ref{Ex: Lorentzprol}, that the  Engel structure $\mathcal{D}$ on $\mathcal{P}C$ is given by
\begin{align*}
    \mathcal{D}_{\Psi(x_1,x_2,x_3,\theta)} & {} = \mathcal{D}_{\tfrac{\cos\theta}{\sqrt{g_{11}}}u_1 + \tfrac{\sin\theta}{\sqrt{g_{22}}}u_2+\tfrac{1}{\sqrt{-g_{33}}}u_3}\\ & {} = (T\pi_L)^{-1}\lspan \frac{\cos\theta}{\sqrt{g_{11}}}u_1 + \frac{\sin\theta}{\sqrt{g_{22}}}u_2+\frac{1}{\sqrt{-g_{33}}}u_3\rspan, 
\end{align*}
	where $\pi_L:\mathcal{P}C\to M$ is the canonical projection. Since $T\pi_L(\partial_{x_i}) = u_i$ and $T\pi_L(\partial_\theta) = 0$, the distribution $\mathcal{D}$ is given by
	\[
	\mathcal{D} = \lspan X:=\frac{\cos\theta}{\sqrt{g_{11}}}\partial_{x_1} + \frac{\sin\theta}{\sqrt{g_{22}}}\partial_{x_2}+\frac{1}{\sqrt{-g_{33}}}\partial_{x_3}, \partial_\theta\rspan.
	\]
	
	Define $\dot{X} : = [\partial_\theta, X] = -\sdfrac{\sin\theta}{\sqrt{g_{11}}}\partial_{x_1}+\sdfrac{\cos\theta}{\sqrt{g_{22}}}\partial_{x_2}$, so the even-contact structure $\mathcal{E}$ on $\mathcal{P}C$ is $\mathcal{E} = \lspan X, \dot{X}, \partial_\theta\rspan$. We have $[X,\dot{X}] = A\partial_{x_1}+B\partial_{x_2}+C\partial_{x_3}$, which is 
	\begin{align*}
	 & \cos^2\theta\Big[\frac{\partial_{x_1}}{\sqrt{g_{11}}}, \frac{\partial_{x_2}}{\sqrt{g_{22}}}\Big]-\sin^2\theta\Big[\frac{\partial_{x_2}}{\sqrt{g_{22}}}, \frac{\partial_{x_1}}{\sqrt{g_{11}}}\Big]
	\\ & -\sin\theta\Big[\frac{\partial_{x_3}}{\sqrt{-g_{33}}}, \frac{\partial_{x_1}}{\sqrt{g_{11}}}\Big] +\cos\theta\Big[\frac{\partial_{x_3}}{\sqrt{-g_{33}}}, \frac{\partial_{x_2}}{\sqrt{g_{22}}}\Big].  
	\end{align*}
	Since $A$ is the $\partial_{x_1}$-component, and analogously for $B$ and $C$, we obtain
	\[
	\begin{cases}
	\begin{aligned}
	A &={} \sdfrac{1}{2g_{11}\sqrt{g_{11}g_{22}}}\sdfrac{\partial g_{11}}{\partial x_{2}}+\sdfrac{\sin\theta}{2 g_{11}\sqrt{-g_{11}g_{33}}}\sdfrac{\partial g_{11}}{\partial x_3},
	\\
	B &={} -\sdfrac{1}{2g_{22}\sqrt{g_{11}g_{22}}}\sdfrac{\partial g_{22}}{\partial x_{1}}-\sdfrac{\cos\theta}{2 g_{22}\sqrt{-g_{22}g_{33}}}\sdfrac{\partial g_{22}}{\partial x_3},\\
	C & ={}-\sdfrac{\sin\theta}{2g_{33}\sqrt{-g_{11}g_{33}}}\sdfrac{\partial g_{33}}{\partial x_1}+\sdfrac{\cos\theta}{2g_{33}\sqrt{-g_{22}g_{33}}}\sdfrac{\partial g_{33}}{\partial x_2}.
	\end{aligned}
	\end{cases}
	\]
	
	Since the kernel $\mathcal{W}$ lies in $\mathcal{D}$ and is not spanned by $\partial_\theta$ (as $[\partial_\theta,\dot{X}]\notin \mathcal{E}$), there exists a smooth function $\mu$ on $\mathcal{P}C$ such that $\mathcal{W} = \lspan X +\mu \partial_\theta\rspan$. We have $[\partial_{\theta}, X +\mu \partial_\theta]\in\mathcal{E}$ and $[X, X +\mu \partial_\theta] \in\mathcal{E}$, whereas
	\[
	[\dot{X}, X +\mu \partial_\theta] = \dot{X}(\mu)\partial_\theta-(A\partial_{x_1}+B\partial_{x_2}+C\partial_{x_3})+\mu\Big(\frac{\cos\theta}{\sqrt{g_{11}}}\partial_{x_1}+\frac{\sin\theta}{\sqrt{g_{22}}}\partial_{x_2}\Big).
	\]

	Since $\dot{X}(\mu)\partial_\theta\in\mathcal{E}$, it is enough to impose that the last two terms belong to $\mathcal{E}$. This is the case if and only if $$-(A\partial_{x_1}+B\partial_{x_2}+C\partial_{x_3})+\mu\Big(\sdfrac{\cos\theta}{\sqrt{g_{11}}}\partial_{x_1}+\sdfrac{\sin\theta}{\sqrt{g_{22}}}\partial_{x_2}\Big)+ C\sqrt{-g_{33}}X\in \mathcal{E}.$$ This vector field reads $$\!-\!\Bigg(\!\!\!\Big(\!A-C\cos\theta\sqrt{\frac{-g_{33}}{g_{11}}}\Big)\partial_{x_1}\!+\Big(\!B-C\sin\theta\sqrt{\frac{-g_{33}}{g_{22}}}\Big)\partial_{x_2}\!\!\Bigg)+\mu\Big(\!\frac{\cos\theta}{\sqrt{g_{11}}}\partial_{x_1}+\frac{\sin\theta}{\sqrt{g_{22}}}\partial_{x_2}\!\Big)\!,$$ and 
	for it to belong to $\mathcal{E}$, it is enough to impose that it is a multiple of $\dot{X}$. We take 
	$$\mu :=A\sqrt{g_{11}}\cos\theta+B\sqrt{g_{22}}\sin\theta-C\sqrt{-g_{33}},$$ so that the vector field becomes $(A\sqrt{g_{11}}\sin\theta-B\sqrt{g_{22}}\cos\theta)\dot{X}$. Hence, by defining $F:=A\sqrt{g_{11}}$, $G:=B\sqrt{g_{22}}$ and $H:=-C\sqrt{-g_{33}}$, the kernel is 
	\begin{align*}
	\mathcal{W} &= \lspan  X+\mu\partial_\theta\rspan  = \lspan Z:=X+(F\cos\theta+G\sin\theta+H)\partial_\theta \rspan.\end{align*}\end{proof}

\begin{Theorem}
\label{Thm: C}
Let $M$ be a three-dimensional spacetime. Then, 
\[
\mathcal{N}  \cong \mathcal{P}C/\mathcal{W}.
\] In addition, if $\mathcal{N}$ is a manifold and $p:\mathcal{P}C\to\mathcal{P}C/\mathcal{W}\cong \mathcal{N}$ is a submersion, the canonical contact structure on $\mathcal{N}$ is
\[
\mathcal{H} \cong p_*\mathcal{E}.
\]
\end{Theorem}

\begin{proof}

We divide the proof of Theorem \ref{Thm: C} into two parts.

\medskip \noindent \textbf{Part I: } $\mathcal{N}  \cong \mathcal{P}C/\mathcal{W}$. It is enough to show that $\mathcal{W}$ computed in Proposition \ref{PropKernel} is pointwise proportional to the geodesic flow $X_g$. Let $(V, \varphi)$ be a local chart of $M$ making the metric diagonal, and let $\gamma(t)= \varphi\big(x_1(t), x_2(t), x_3(t)\big)$ be a null geodesic. Then,
\begin{equation}
\label{eq: conditionDerivatives}
(x_1')^2g_{11}+(x_2')^2g_{22} = -(x_3')^2g_{33}
\end{equation}
and the curve $x_1'u_1+x_2'u_2+x_3'u_3\in T_{\varphi(x_1, x_2, x_3)}M$ is an integral line of $X_g$. Under our identification $\mathcal{P}C\cong\{\sdfrac{\cos\theta}{\sqrt{g_{11}}}u_1 + \sdfrac{\sin\theta}{\sqrt{g_{22}}}u_2+\sdfrac{1}{\sqrt{-g_{33}}}u_3\}_{\theta\in\S^1, x\in M}$, this curve reads
\[
\frac{x_1'}{x_3'}\sqrt{\frac{g_{11}}{-g_{33}}}\frac{u_1}{\sqrt{g_{11}}}+\frac{x_2'}{x_3'}\sqrt{\frac{g_{22}}{-g_{33}}}\frac{u_2}{\sqrt{g_{22}}}+\frac{u_3}{\sqrt{-g_{33}}}\in T_{\varphi(x_1, x_2, x_3)}M = \Psi\big(x_1, x_2, x_3, \theta\big),
\]
where $\theta$ is such that $\cos\theta = \frac{x_1'}{x_3'}\sqrt{\frac{g_{11}}{-g_{33}}}$ and $\sin\theta = \frac{x_2'}{x_3'}\sqrt{\frac{g_{22}}{-g_{33}}}$. Note that this makes sense because of Equation \eqref{eq: conditionDerivatives}, and that $x_3'$ is nonzero if $\gamma$ is nonconstant. Then, the tangent vector to this curve, which gives the expression of $X_g$ in the coordinate chart $\big(V\times(0, 2\uppi)\big)$ on $\mathcal{P}C$, is
\begin{align*}
x_1'\partial_{x_1}+&x_2'\partial_{x_2}+x_3'\partial_{x_3}+\theta'\partial_{\theta} \\ &=  x_3'\sqrt{-g_{33}}\Big(\frac{\cos\theta}{\sqrt{g_{11}}}\partial_{x_1}+\frac{\sin\theta}{\sqrt{g_{22}}}\partial_{x_2}+\frac{1}{\sqrt{-g_{33}}}\partial_{x_3}+\frac{\theta'}{x_3'\sqrt{-g_{33}}}\partial_\theta\Big).
\end{align*}

Hence, our claim is equivalent to $\frac{\theta'}{x_3'\sqrt{-g_{33}}} = F\cos\theta+G\sin\theta + H$. Using that $-\theta'\sin\theta = \Big(\frac{x_1'}{x_3'}\sqrt{\frac{g_{11}}{-g_{33}}}\Big)'$ and that the $x_i$ satisfy the geodesic equation the result follows, after computing the Christoffel symbols for the metric in Equation \eqref{matrixdiag}.

\medskip
\noindent \textbf{Part II: } $\mathcal{H} \cong p_*\mathcal{E}$. Firstly, as argued in \cite[p. 246]{Adachi_2002}, the even-contact structure $\mathcal{E}$ is invariant under the flow of any vector field generating $\mathcal{W}$. Therefore, the pushforward $p_*\mathcal{E}$ is well defined.

	Let $\gamma\in\mathcal{N} \cong \mathcal{P}C/ \mathcal{W}$. Then, $\gamma$ is given by a curve $\mu: (-\varepsilon, \varepsilon)\to M$ which is a null geodesic in $M$, and $\gamma = p\big(\dot{\mu}(s)\big)$ for all $s\in(-\varepsilon, \varepsilon)$. Let $q_0 := \mu(0)$ and define a coordinate system $\varphi:(x_1,x_2,x_3)\mapsto\varphi(x_1,x_2,x_3)$ around $q_0$ given by Proposition \ref{Prop: DiagonalMetric}. Let $\Psi(x_1,x_2,x_3, \theta)$ be coordinates around $\dot{\mu}(0)$ in $\mathcal{P}C$ as defined before Proposition \ref{PropKernel}. Now, for all $s\in[0,\varepsilon)$ small enough, the point $q_s: = \mu(s)$ lies in the image of $\varphi$ and $\dot{\mu}(s)$ lies in the image of $\Psi$, and we can consider functions $x_i,\theta$ such that $\dot{\mu}(s) = \Psi\big(x_1(s), x_2(s), x_3(s), \theta(s)\big)$. If $s>0$ is small enough, the points $q_0$ and $q_s$ are not conjugate in $M$.
	
	The sky of $q_0$ is given by
	\[
	\mathfrak{S}_{q_0} = \{p\circ\Psi\big(x_1(0), x_2(0), x_3(0), \theta\big)\ |\ \theta\in \big(\theta(0)-\uppi, \theta(0)+\uppi\big]\},
	\]
	and similarly for $\mathfrak{S}_{q_s}$. Therefore, $$\mathcal{H}_\gamma= T_\gamma\mathfrak{S}_{q_0}\oplus T_\gamma\mathfrak{S}_{q_s}= T_{\dot{\mu}(0)}p(\lspan \partial_\theta\rspan)\oplus T_{\dot{\mu}(s)}p(\lspan \partial_\theta\rspan).$$ We ought to express the second addend as a pushforward of a line over $\dot{\mu}(0)$. Take $s>0$ small enough and let $K$ be a neighbourhood of $\dot{\mu}(s)$ for which the flow $\Phi_{-s}^Z:\nolinebreak K \to  \Phi_{-s}^Z(K)$ at time $-s$ of the vector field $Z =  X + (F\cos\theta+G\sin\theta)\partial_\theta$ is a diffeomorphism. Then, $p\circ\Phi^Z_{-s} = p$, and since $\Phi^Z_{-s}\big(\dot{\mu}(s)\big) = \dot{\mu}(0)$, we can compute
	\[
	T_\gamma\mathfrak{S}_{q_s} = T_{\dot{\mu}(s)}p(\partial_\theta) = (T_{\dot{\mu}(0)}p\circ T_{\dot{\mu}(s)}\Phi^{Z}_{-s})(\partial_\theta)
	\]
	and
	\begin{align*}
		\mathcal{H}_{\gamma} &= T_{\dot{\mu}(0)}p\left(\lspan \partial_\theta,  T_{\dot{\mu}(s)}\Phi^{Z}_{-s}(\partial_\theta) \rspan\right) = T_{\dot{\mu}(0)}p\left(\lspan \partial_\theta,  T_{\dot{\mu}(s)}\Phi^{Z}_{-s}(\partial_\theta)-\partial_\theta \rspan\right) \\ &= T_{\dot{\mu}(0)}p\left(\lspan \partial_\theta, \frac{ T_{\dot{\mu}(s)}\Phi^{Z}_{-s}(\partial_\theta)-\partial_\theta}{s} \rspan\right),
	\end{align*}
	for all $s>0$ small enough. Hence, the result is still true if we take the limit $s\to 0$. Thus, we obtain $\mathcal{H}_{\gamma}  = T_{\dot{\mu}(0)}p\left(\lspan \partial_\theta, \lim\limits_{s\to 0}\frac{ T_{\dot{\mu}(s)}\Phi^{Z}_{-s}(\partial_\theta)-\partial_\theta}{s} \rspan\right) =  T_{\dot{\mu}(0)}p\left(\lspan\partial_\theta,  [\partial_\theta, Z]\rspan\right)$. We compute $[\partial_\theta, Z] = \dot{X}+(-F\sin\theta+G\cos\theta )\partial_\theta$ and therefore
	\begin{align*}
		\mathcal{H}_{\gamma}  &= T_{\dot{\mu}(0)}p\left(\lspan\partial_\theta,  [\partial_\theta, Z]\rspan\right) = T_{\dot{\mu}(0)}p\left(\lspan\partial_\theta,  \dot{X} \rspan\right) = T_{\dot{\mu}(0)}p\left(\lspan\partial_\theta,  \dot{X}, Z \rspan\right) \\&  = T_{\dot{\mu}(0)}p\left(\lspan\partial_\theta,  \dot{X}, X+ (F\cos\theta+G\sin\theta)\partial_\theta \rspan\right) \\ &=T_{\dot{\mu}(0)}p\left(\lspan\partial_\theta,  \dot{X}, X \rspan\right) = T_{\dot{\mu}(0)}p(\mathcal{E}).
	\end{align*}
	
	This concludes the proof that $\mathcal{H} \cong p_*\mathcal{E}$ and hence the theorem.\end{proof}
	
\begin{Example}\label{ex: S2timesS1Deprolong}
The Lorentzian prolongation of the spacetime $(\S^2\times\S^1,g_1)$ is diffeomorphic to $\mathcal{P}C\cong ST\S^2\times\S^1$. Letting $t$ be the $\S^1$-coordinate and $\theta$ the coordinate on the fibre of $ST\S^2$, the Engel flag over $(u,s)\in ST\S^2\times\S^1$ is, seeing $T\S^2\subset T(ST\S^2)$ via the Levi-Civita connection on $\S^2$,
\[
\lspan u+\partial_t\rspan \subset \lspan u+\partial_t, \partial_\theta\rspan \subset  \lspan u\rspan^\perp\oplus\lspan u+\partial_t,\partial_\theta\rspan\subset T_{(u,s)}(ST\S^2\times\S^1).
\]
\noindent Hence, $(ST\S^2\times\S^1)/\mathcal{W}  = (ST\S^2\times\S^1)/\lspan u+\partial_t\rspan \cong ST\S^2\cong\mathcal{N}_1$, and letting $p: ST\S^2\times\S^1\to ST\S^2$ be the natural projection, $p_*\mathcal{E}\cong \chi\cong\mathcal{H}$.
\end{Example}


\subsection{Retrieving the spacetime}\label{sec:recovering}
We investigate how Theorem \ref{Thm: C} allows us to characterize the three-dimensional contact manifolds that are spaces of null geodesics of a spacetime. This provides a procedure to, given a contact manifold $(N, \xi)$ satisfying the necessary conditions, find a spacetime whose space of null geodesics is precisely $(N, \xi)$.

The main idea is to  Cartan-prolong the contact manifold and make use of Remark~\ref{Rk: LorentzExample} to Lorentz-deprolong it. 
\begin{equation}\label{diag:Engelx2}
\begin{tikzcd}[column sep=small]
& (\mathcal{P}C,\Dtilde )\arrow[drrr,swap, "\text{Thm.} \ref{Thm: C}"] \arrow[rr,leftrightarrow, dashed, "?"]
& & (\S(\xi),\mathcal{D})\arrow[dlll, dashed, "?"]& \\
(M,g) \arrow[ur, "\text{Ex.}\ref{Ex: Lorentzprol}"] \arrow[rrrr, swap,"\text{Def.}\ref{Def: SpaceofNullGeod}"] & & & & (N,\xi) \arrow[ul,swap, ,"\text{Ex.}\ref{Ex: CartanProl}"]
\end{tikzcd}
\end{equation}

We have the following result.

\begin{Theorem}
\label{thm: characterisationEngelFoliation}
A three-dimensional contact manifold $(N, \xi)$ is contactomorphic to the space of null geodesics of a spacetime if and only if there exists an Engel manifold $(Q,\mathcal{D})$ with Engel flag $\mathcal{W}\subset\mathcal{D}\subset\mathcal{E}\subset TQ$ such that
\begin{equation}
\label{eq: ConditionQuotient}
N\cong Q/\mathcal{W} \hspace{1cm}\text{and}\hspace{1cm}\xi\cong p_*\mathcal{E},
\end{equation}
for $p:Q\to Q/\mathcal{W}$ the projection, and $Q$ admits an oriented foliation by circles $\mathcal{F}$ such that
\begin{enumerate}
    \item for all $S\in\mathcal{F}$ and $x\in S$, we have $T_xS\oplus \mathcal{W}_x = \mathcal{D}_x$,
    \item the space of leaves $M := Q/\mathcal{F}$ is a manifold and the projection $q:~Q\to M$ is a submersion,
    \item for every $S\in\mathcal{F}$, the image $q_*\mathcal{D}|_S$ is a cone in the vector space $T_{q(S)}M$ and the map $x\in S\mapsto q_*\mathcal{D}_x$ is injective.
\end{enumerate}

In addition, if the above conditions are satisfied, $(N, \xi)$ is contactomorphic to the space of null geodesics of $(Q/\mathcal{F}, g)$, where $g$ is a metric on $Q/\mathcal{F}$ with bundle of cones $q_*\mathcal{D}$.
\end{Theorem}
\begin{proof}
Let $(N, \xi)$ be the space of null geodesics of a spacetime $(L, h)$. The existence of $Q$ satisfying \eqref{eq: ConditionQuotient} follows from Theorem \ref{Thm: C} taking $(Q,\mathcal{D})$ as the Lorentz prolongation of $(L, h)$. Then, we can take $\mathcal{F} := \{\mathcal{P}C_x\}_{x\in L}$, which is a foliation by circles and is oriented because $L$ is a spacetime. By definition of the Lorentz prolongation, $T\mathcal{F}\subset \mathcal{D}$, and by Example \ref{Ex: Lorentzprol}, the kernel $\mathcal{W}$ is transverse to $T\mathcal{F}$. Hence $i)$ is satisfied. Also, by Remark \ref{Rk: LorentzExample}, $Q/\mathcal{F}\cong L$, and therefore it is a manifold, and $q$ is the projection $\pi_L:Q=\mathcal{P}C\to L$, which is a submersion. Hence, $ii)$ follows.

Also by definition of the Lorentz prolongation, for all $u\in Q$, $q_*\mathcal{D}_u = u$, where on the right hand side we regard $u$ as a vector on a cone of $L$. Therefore, $q_*\mathcal{D}|_{\mathcal{P}C_x} = C_x$ and $iii)$ also follows. Finally, by Theorem \ref{Thm: C} and \eqref{eq: ConditionQuotient}, the last claim is also satisfied, as any other metric on $L$ with the same bundle of cones is conformal to $g$, and therefore has the same space of null geodesics.

For the converse, let $(N, \xi)$ be a contact manifold such that $N = Q/\mathcal{W}$ and $\xi = p_*\mathcal{E}$ for an Engel manifold $Q$ with flag $\mathcal{W}\subset\mathcal{D}\subset\mathcal{E}\subset TQ$. Assume $Q$ admits an oriented foliation $\mathcal{F}$ satisfying $i),\ ii),\ iii)$ above. Let $M:=Q/\mathcal{F}$, which is a manifold by hypothesis, and $q:Q\to M$ the projection map. Now, since $T\mathcal{F}\subset\mathcal{D}$, the pushforward $ q_* \mathcal{D}_x$ is a line in $T_{q(x)}M$ for all $u\in \mathcal{P}C$ and by $iii)$, these create a cone in $T_{q(S)}M$ when traveling the leaf $S\in\mathcal{F}$ containing $x$. Therefore, we obtain a smooth bundle of cones on $TM$, and hence there exists a metric $g$ on $M$ with such bundle of cones. In addition, any two metrics with the same bundle of null cones are conformal, and hence produce the same contact manifold of null geodesics. Since $\mathcal{F}$ is oriented, we can assign a consistent orientation to each cone in the bundle and therefore $(M,g)$ is a spacetime.

Let now $(\mathcal{P}C, \Dtilde)$ be the Lorentz prolongation of $(M, g)$, with flag $\Wtilde\ \subset\ \Dtilde\ \subset\ \Etilde\ \subset T(\mathcal{P}C)$, and define the projection $\pi_M:\mathcal{P}C\to M$. Consider the map
\[ \begin{array}{cccc}
         \Phi:& Q&\to &\mathcal{P}C   \\
         &x  &\mapsto  &q_*\mathcal{D}_{x}
    \end{array},
   \]
  which is a diffeomorphism by the definition of $g$ and the hypothesis that $x\in S\mapsto p_*\mathcal{D}_x$ is injective for any $S\in\mathcal{F}$. In addition, 
  \[
\Big({\pi_L}_*\circ\Phi_*\Big)(\mathcal{D}_{x}) = (\pi_M\circ\Phi)_*(\mathcal{D}_{x}) = q_*\mathcal{D}_{x} = \Phi(x),
\]
which implies $\Phi_*\mathcal{D}\subset\ \Dtilde$. Since $\Phi$ is a submersion, we find that $\Phi_*\mathcal{D} = \Dtilde$, and so $\Phi$ is an Engel-morphism. Let us define $r:\mathcal{P}C\to\mathcal{P}C/\Wtilde$ the projection. Then, by Theorem \ref{Thm: C}, we find that the space of null geodesics $(\mathcal{N}, \mathcal{H})$ of $(M, g)$ is
\[
\mathcal{N}\cong\mathcal{P}C/\Wtilde\ \cong Q/\mathcal{W}\cong M 
\]
and
\[
\mathcal{H}\cong r_*\Etilde\ \cong p_*\mathcal{E}\cong \xi, 
\]
where we make use of Remark \ref{Rk: CartanExample}.
\end{proof} 

\begin{Remark}
We have later found out that \cite[Rk. 1.7]{Mitsumatsu} suggests the necessity of a foliation $\mathcal{F}$ like the one in Theorem \ref{thm: characterisationEngelFoliation}. 
\end{Remark}

Since the Cartan prolongation of $(N,\xi)$ satisfies \eqref{eq: ConditionQuotient}, we obtain:

\begin{Corollary}
\label{thm: retrieveNEW}
A three-dimensional contact manifold $(N, \xi)$ is contactomorphic to the contact manifold of null geodesics of a spacetime if the Cartan prolongation $(\S(\xi), \mathcal{D})$ of $(N, \xi)$ admits an oriented foliation $\mathcal{F}$ by circles such that
\begin{enumerate}
    \item for all $S\in\mathcal{F}$ and $u\in S$, we have $T_uS\oplus \mathcal{W}_u = \mathcal{D}_u$, where $\mathcal{W}$ denotes the kernel of $(\S(\xi), \mathcal{D})$,
    \item the space of leaves $M := \S(\xi)/\mathcal{F}$ is a manifold and the projection $p:~\S(\xi)\to M$ is a submersion,
    \item for every $S\in\mathcal{F}$, the image $p_*\mathcal{D}|_S$ is a cone in the vector space $T_{p(S)}M$ and the map $u\in S\mapsto p_*\mathcal{D}_u$ is injective.
\end{enumerate}
In addition, if $i),\ ii),\ iii)$ are satisfied, $(N, \xi)$ is contactomorphic to the space of null geodesics of $(M,g)$, where $g$ is any metric on $M$ with bundle of null cones $p_*\mathcal{D}$.
\end{Corollary}


We continue our discussion by exploring the relation between the Cartan prolongation of $(N,\xi)$ and the Engel manifold $Q$ in Theorem \ref{thm: characterisationEngelFoliation}, whenever it exists. The following proposition is an adaptation of \cite[Prop. 5.4]{Montgomery_2001}.
 
\begin{Proposition}
\label{prop: localdiffeo}
Let $(Q,\mathcal{D})$ be an Engel manifold with flag $\mathcal{W}\subset\mathcal{D}\subset\mathcal{E}\subset TQ$. Assume $Q/\mathcal{W}$ is a manifold, which then can be endowed with a contact structure $\xi:=p_*\mathcal{E}$ for $p:Q\to Q/\mathcal{W}$ the projection, provided it is a submersion. Then, there exists a local diffeomorphism $\Phi:Q\to \S(\xi)$ to the Cartan prolongation which is compatible with the Engel structure.
\end{Proposition}
\begin{proof}
Let $x\in X$. Since $\mathcal{W}_x\subset\mathcal{D}_x$, the pushforward $p_* \mathcal{D}_x$ is a line in $\xi_{p(x)}$. Hence, the map
\[
\begin{array}{cccc}
     \Phi:&Q&\to&\S(\xi)  \\
     & x&\mapsto& p_*\mathcal{D}_x
\end{array}
\]
is well defined, and smooth. Let $x\in Q$ and consider a neighbourhood $U\subset Q$ of $x$ in which we can trivialize $\mathcal{D}|_U = \lspan Z, Y\rspan$, for $Z\in\mathfrak{X}(U)$ spanning $\mathcal{W}$. Since $\Phi$ is a bundle map over $Q/\mathcal{W}$, it is enough to show that $T_y\Phi$ is surjective when restricted to $\mathcal{W}_y$ for every $y\in U$. Therefore, by linearity, it is enough to show $\Phi_*Z\neq 0$. By definition of $\Phi$, this is equivalent to $[Z,Y]\neq 0$, which holds because $\mathcal{D}$ is Engel. 

It is only left to show that $\Phi$ preserves the Engel structure. Let $u\in\mathcal{D}_x$ and define $\pi_C:\S(\xi)\to X/\mathcal{W}$ the projection. Then,
\[
\Big({\pi_C}_*\circ \Phi_*\Big)(\mathcal{D}_x) = (\pi_C\circ \Phi)_*(\mathcal{D}_x) = p_*\mathcal{D} = \Phi(x).
\]
Since $\Phi$ is a submersion, the Proposition follows.
\end{proof}

\begin{Corollary}
\label{cor: LocalDiffeo}
If $(N, \xi)$ is a three-dimensional contact manifold contactomorphic to the space of null geodesics of a spacetime, the Engel manifold $(Q,\mathcal{D})$ described in Theorem \ref{thm: retrieveNEW} comes with a canonical local diffeomorphism $\Phi: Q\to\S(\xi)$ compatible with the Engel structures.
\end{Corollary}

Diagram \eqref{diag:Engelx2} above is thus completed to 
\begin{equation*}
\begin{tikzcd}[column sep=small]
& (\mathcal{P}C,\Dtilde )\arrow[drrr, swap,"\text{Thm.} \ref{Thm: C}", near end] \arrow[rr,rightarrow, "\text{Cor.}\ref{cor: LocalDiffeo}"]
&  & (\S(\xi),\mathcal{D})\arrow[dlll,"\text{Thm.} \ref{thm: characterisationEngelFoliation}"]& \\
(M,g) \arrow[ur, "\text{Ex. 4.4}"] \arrow[rrrr, swap,"\text{Space of null geodesics}"] & & & & (N,\xi). \arrow[ul,swap, ,"\text{Ex. 4.2}"]
\end{tikzcd}
\end{equation*}
We believe that this approach can be useful in order to answer the open question of whether the contact structure of the space of null of geodesics can be overtwisted, but describing or just dealing with the foliation in Theorem \ref{thm: characterisationEngelFoliation} will require further work. 

We finally look at two illustrative examples where the manifolds involved can be described explicitly and the subtleties of the main results can be appreciated.

\begin{Example}\label{ex:R3S1}
Consider $\R^3$ with coordinates $(x, y, t)$ and tangent vector fields $u_x, u_y, u_t$. Then, the Lorentz prolongation of $(\R^3, dx^2+dy^2-dt^2)$ is diffeomorphic to $\R^3\times\S^1$ via $(x,y,t,\theta)\mapsto \lspan\cos\theta u_x+\sin\theta u_y+u_t\rspan\subseteq T_{(x,y,t)}\R^3$. 
Under this identification, the Engel structure is $\mathcal{D}_{(x,y,t,\theta)} = \lspan \cos\theta\partial_x+\sin\theta\partial_y+\partial_t, \partial_\theta\rspan$, which implies $\mathcal{W} = \lspan \cos\theta\partial_x+\sin\theta\partial_y+\partial_t\rspan$ and $\mathcal{E} = \mathcal{D}\oplus\lspan-\sin\theta\partial_x+\cos\theta\partial_y\rspan$. The foliation $\mathcal{F}$ by circles is given by $\{\{k\}\times\S^1\}_{k\in\R^3}$. We also have $\R^3\times\S^1/\mathcal{W}\cong \R^2\times\S^1$ with projection map
\[
\begin{array}{cccc}
     p:&\R^3\times\S^1&\to & \R^2\times\S^1 \\
     &(x,y,t,\theta)&\mapsto& (x-t\cos\theta, y-t\sin\theta, \theta). 
\end{array}
\]

Take coordinates $(u,v,\theta)$ on $\R^2\times\S^1$ with tangent vectors $w_u, w_v, w_\theta$. Then, the induced contact structure on $\R^2\times\S^1$ is
\begin{align*}
\xi_{p(x,y,t,\theta)} :=p_*\mathcal{E}_{(x,y,t,\theta)} &= \lspan -\sin\theta w_u+\cos\theta w_v, t\sin\theta w_u-t\cos\theta w_v+ w_\theta\rspan \\ &= \lspan -\sin\theta w_u+\cos\theta w_v, w_\theta\rspan, 
\end{align*}
which is the canonical contact structure on $\R^2\times\S^1\cong ST\R^2$. The Cartan prolongation $\S(\xi)$ of $(\R^2\times\S^1, \xi)$ is diffeomorphic to $\R^2\times\S^1\times\S^1$ via $$(u,v,\theta, \omega)\mapsto \orlspan \cos\omega(-\sin\theta w_u+\cos\theta w_v)+\sin\omega w_\theta\orrspan\in \S(\xi)_{(u,v,\theta)}.$$ Under this identification, the Engel distribution reads $$\Dtilde_{(u,v,\theta, \omega)} = \lspan \cos\omega(-\sin\theta \partial_u+\cos\theta \partial_v)+\sin\omega \partial_\theta, \partial_{\omega}\rspan \subset T_{(u,v,\theta,w)}\S(\xi).$$

The local diffeomorphism $\Phi$ defined in Proposition \ref{prop: localdiffeo} is $(x,y,t,\theta)\mapsto p_*\partial_\theta = \orlspan t\sin\theta w_u-t\cos\theta w_v+w_\theta\orrspan = \orlspan \frac{t}{\sqrt{1+t^2}}\sin\theta w_u-\frac{t}{\sqrt{1+t^2}}\cos\theta w_v+\frac{1}{\sqrt{1+t^2}}w_\theta\orrspan$. Hence, using the identification above,
\[
\begin{array}{cccc}
     \Phi:&\R^3\times\S^1&\to&\R^2\times\S^1\times\S^1  \\
     &(x,y,t,\theta)&\mapsto & (x-t\cos\theta, y-t\sin\theta, \theta,\arccos{\frac{t}{\sqrt{1+t^2}}} )
\end{array},
\]
where we take $\arccos\frac{t}{\sqrt{1+t^2}}\in(0, \uppi)$. Then, $\Phi$ is a global diffeomorphism onto its image $\Phi(\R^3\times\S^1) = \R^2\times\S^1\times(0,\uppi)$. Note that the image of $\mathcal{F}$ under $\Phi$ is $\{\Phi(x,y,t, \theta)\ |\ \theta\in\S^1\}_{(x,y,t)\in\R^3}$. This is still a foliation by circles on the image of $\Phi$ satisfying $i)$ and $ii)$ above with projection
\[
\begin{array}{cccc}
     q:&\R^2\times\S^1\times(0,\uppi)&\to&\R^2\times\S^1\times(0,\uppi)/\Phi(\mathcal{F})\cong \R^2\times(0,\uppi)  \\
     &(u,v,\theta,\omega)&\mapsto&(u+\cos\theta\frac{\cos\omega}{\sin\omega}, v+\sin\theta\frac{\cos\omega}{\sin\omega}, \omega).
\end{array}
\]

Again, taking $(u,v,\omega)$ as coordinates on $\R^2\times(0,\uppi)$, we can compute
\[
q_*\Dtilde = q_*\partial_\omega = \lspan \partial_\omega-\frac{\cos\theta}{\sin^2\omega}\partial_u  -\frac{\sin\theta}{\sin^2\omega}\partial_v\rspan =: V(\omega, \theta),
\]
which indeed defines a bundle of cones with injective map $\theta\in\S^1\mapsto V(\omega,\theta)$ for a given $\omega\in(0,\uppi)$. Hence, $iii)$ is also satisfied.

The bundle of cones obtained on $\R^2\times(0,\uppi)$ is the induced by the metric $g = du^2+dv^2-\sin^4\omega d\omega^2$, and therefore we obtain a spacetime isometric to $(\R^3, dx^2+dy^2-dt^2)$ via $(x,y,t)\mapsto(x,y,\arccos{\frac{t}{\sqrt{1+t^2}}})$.
\end{Example}

\begin{Example}
Let $(M, g) = (\S^2\times\S^1, g_\circ-dt^2)$. As argued in Example \ref{ex: S2timesS1Deprolong}, $\mathcal{P}C\cong ST\S^2$ with flag 
\[
\lspan u+\partial_t\rspan \subset \lspan u+\partial_t, \partial_\theta\rspan \subset  \lspan u\rspan^\perp\oplus\lspan u+\partial_t,\partial_\theta\rspan\subset T_{(u,s)}(ST\S^2\times\S^1),
\]
and $(\mathcal{N}, \mathcal{H}) \cong (ST\S^2, \chi)$. Note that, since the $\S^1$-bundle $ST\S^2\to \S^2$ is oriented, the vector field $\partial_\theta$ is well defined everywhere. In addition, for $u\in ST_x\S^2$, we can define $u^\perp\in ST_x\S^2$ orthogonal to $u$ with respect to $g_\circ$ and so that $(u, u^\perp, x)$ is a positive basis for $\S^2$. Seeing $T\S^2\subset T(ST\S^2)$ via the Levi-Civita connection on $(\S^2,g_\circ)$, we can see $u^\perp\in T_u(ST\S^2)$, and we obtain a well-defined vector field $P$ on $ST\S^2$ given pointwise by $P_u = u^\perp$. Then, the contact structure on $ST\S^2$ is  $\chi = \lspan\partial_\theta, P\rspan$ and, in particular, it is trivial as a vector bundle. Therefore the Cartan prolongation $\S(\chi)$ of $(ST\S^2, \chi)$ is $\S(\chi)\cong ST\S^2\times\S^1$ via $(u,s)\mapsto \orlspan \cos{s}P + \sin{s}\partial_\theta\orrspan\in T_u ST\S^2$. 

Now, the local diffeomorphism $\Phi:ST\S^2\times\S^1\to ST\S^2\times\S^1$ is given by
\[
\Phi(u,s) = p_*\mathcal{D}_{u,s} = p_*\partial_\theta,
\]
which can be described as follows. The element in $ST\S^2$ of $\Phi(u,s)$ is simply $p(u,s)$, which recall can be described as taking the great circle $\mu\subset\S^2$ defined by $u$, and parallel transporting $u$ over $\mu$ an angle $\theta$ backwards. We get
\[
\Phi(u, s) = \orlspan \cos{s}\partial_{\theta}-\sin{s}\ p(u,s)^\perp \orrspan\in T_{p(u,s)}ST\S^2,
\]
that is, $\Phi(u, s) = \big(p(u, s), s+\frac{\pi}{2}\big)$. This is a global diffeomorphism $\Phi:~ST\S^2\to ST\S^2$ preserving the Engel structures. The foliation $$\mathcal{F} = \{ST_x\S^1\times\{t\}\}_{(x,t)\in \S^2\times\S^1}$$ on $\mathcal{P}C\cong ST\S^2\times\S^1$ gets sent under $\Phi$ to a foliation by circles $\mathcal{F}_\Phi$ which can be described as follows. The foliation splits in a family of $\S^1-$foliations of $ST\S^2$. On $ST\S^2\times\{t\}$, $\mathcal{F}_\Phi$ is formed by all circles on $\S^2$ with radius $t$ and vectors pointing towards the centre of such circle. Hence, $\mathcal{F}_\Phi$ satisfies $i)$. Now, fixing $t$, every circle of $\mathcal{F}_\Phi$ on $ST\S^2\times\{t\}$ has a unique centre, and every point of $\S^2$ defines one such circle. Therefore, $ST\S^2/\mathcal{F}_\Phi\cong \S^2\times\S^1$, and the projection $q:ST\S^2\to \S^2\times\S^1$ is a submersion, since it is on every slice $q|_{ST\S^2\times\{t\}}:ST\S^2\times\{t\}\to \S^2\times\{t\}$, as this is only parallel transporting a tangent vector an angle $t$ over its great circle, and taking the basepoint. Therefore, $ii)$ is also satisfied. The bundle of cones defined on $\S^2\times \{t\}$ is that of a strictly positive constant multiple of $g_c$, and this constant varies smoothly with respect to $t$. Hence, we obtain a spacetime isometric to $(\S^2\times\S^1, g_c\times dt^2)$.

\end{Example}


\end{document}